\documentclass[11pt]{article}
\usepackage{amsmath,amssymb}
\usepackage{graphicx}
\usepackage{epsfig}
\usepackage{color}
\usepackage{float}
\usepackage{ulem}
\newcommand{\cqfd}%
{\mbox{}\nolinebreak\hfill\rule{2mm}{2mm}\medskip\par}

\newcommand{\prodscal}[2]{\left( \, {#1} \, | \, {#2} \, \right)}

\newcommand{\norml}[2]{\left\| {#1} \right\|_{L^2({#2})}}

\newcommand{\derp}[2]{\frac{\partial{#1}}{\partial{#2}}}

\def\divg{\mathop{\rm div}\nolimits}
\def\cqfd{\qed}
\def\IR{\mathbb{R}}

\def\vtld{\widetilde{v}}

\def\nabtld{\widetilde{\nabla}}

\def\divtld{\widetilde{\divg}}

\newtheorem{theo}{Theorem}[section]
\newtheorem{lemm}[theo]{Lemma}

\newtheorem{rem}[theo]{Remark}

\newenvironment{dem}[1] {\par\noindent{\it Proof. }{#1}}{$\square$}
\begin{document}
\title{A Transport for imaging process}
\author{Olivier Besson \\
Universit\'e de Neuch\^atel, \\
Institut de Math\'ematiques \\
11, rue E. Argand, 2000 Neuch\^atel, Switzerland, \\
olivier.besson@unine.ch \smallskip \\
Martine Picq \\
Universit\'e de Lyon CNRS  \\
INSA-Lyon ICJ UMR 5208, bat. L. de Vinci, \\
20 Av. A. Einstein, F-69100 Villeurbanne Cedex France \\
martine.picq@insa-lyon.fr \smallskip \\
J\'er\^ome Pousin \\
Universit\'e de Lyon CNRS  \\
INSA-Lyon ICJ UMR 5208, bat. L. de Vinci, \\
Universit\'e de Lyon CNRS  \\
20 Av. A. Einstein, F-69100 Villeurbanne Cedex France \\
jerome.pousin@insa-lyon.fr
\date{June 2012}}
\maketitle
\begin{abstract}
This work originates from a heart's images tracking which is to
generate an apparent continuous motion, observable through intensity
variation from one starting image to an ending one both supposed
segmented. Given two images $\rho_0$ and $\rho_1$, we calculate an
evolution process $\rho(t,\cdot)$ which transports $\rho_0$ to
$\rho_1$ by using the optical flow. In this paper we propose
an algorithm based on a fixed point formulation and a space-time least squares
formulation of the transport  equation for computing a
transport problem. Existence results are given for a transport problem with a minimum 
divergence for a dual norm or a weighted $H^1_0$-semi norm, for the velocity. 
The proposed transport is compare with the transport introduced by Dacorogna-Moser. The strategy is
implemented in a 2D case and numerical results are presented with a first order Lagrange
finite element, showing the efficiency of the proposed strategy.
\end{abstract}
\textbf{keywords} AMS Classification 35F40; 35L85; 35R05;
\section{Introduction}
\label{Intro}
Modern medical imaging modalities can
provide a great amount of information to study the human anatomy and
physiological functions in both space and time. In cardiac magnetic
resonance imaging (MRI) for example, several slices can be acquired
to cover the heart in 3D and at a collection of discrete time
samples over the cardiac cycle.  From these partial observations, the challenge is to
extract the heart's dynamics from these input spatio-temporal data throughout the cardiac cycle
\cite{Lynch08}, \cite{Schaerer08}.

Image registration consists in estimating a transformation which insures the warping
of one reference image onto another target image (supposed to present some similarity).
Continuous transformations are privileged, the sequence of transformations during the
estimation process is usually not much considered. Most important is the final resulting
 transformation and not the way one image will be transformed to the other.
Here, we consider a reasonable interpolation process to continuously map the image intensity
functions between two
images in the context of cardiac motion estimation and modeling.

The aim of this paper is to present, in the context of
optical flow, an algorithm to compute a time dependent transportation
plan without using lagrangian techniques.

The paper is organized as follows. The introduction is ended, by recalling
the optical flow model (OF) . In section \ref{algo}, the algorithm is presented, and
its convergence is discussed.
In section \ref{masse} it is shown that the solutions obtained with the
proposed algorithm are the solutions minimizing the same energy than the
time dependent optimal mass transportation problem.  Section \ref{exemples} is devoted to
numerical results. In particular a 2D cardiac  medical image  is considered.

\subsection{The optical flow (OF) method}
Let $\rho$ be  the intensity function, and  $v$ be the
velocity of the apparent motion of brightness pattern. An image
sequence is considered via the gray-value map $\rho : Q=(0,1) \times
\Omega \rightarrow \IR$  where $\Omega \subset \IR^d$ is a bounded
regular domain, the support of images, for $d=1,2,3$. If image
points move according to the velocity field $v: \, Q \rightarrow
\IR^d$, then gray values $\rho(t,X(t,x))$ are constant along motion
trajectories $X(t,x)$. One obtains the optical flow equation.
\begin{equation}
 \frac{d}{dt} \rho(t,X(t,x))= \partial_t \rho(t,X(t,x))  +
 \prodscal{v}{\nabla_X \rho(t,X(t,x))}_{\IR^d}  = 0.
\end{equation}
The previous equations lead to an ill-posed problem for the unknown
$(\rho,v)$. Variational formulations or relaxed minimizing problems
for computing jointly $(\rho,v)$ have been first proposed in
\cite{Aubert:1999} and after by many other authors. Here our concern
is somewhat  different. Finding $(\rho,v)$ simultaneously  is
possible by solving a mass transport problem. Similarly to the work  developed in
\cite{Benamou-Brenier,Brenier}, a characterization of  $(\rho,v)$ as solution of a minimizing
problem is  developed.

Let $\rho_0$ and $\rho_1$ be the cardiac images
between two times arbitrary fixed to zero and one, the mathematical
problem reads: find $\rho$ the gray level function defined from $Q$
with values in $[0,1]$ verifying
\begin{equation}\label{optiproblem}
\left \{ \begin{array}{l}
\partial_t \rho(t,x)  + \prodscal{(v(t,x)}{\nabla \rho(t,x)} = 0, \, (t,x) \in
(0,1)\times \Omega\\
\rho(0,x)=\rho_0(x); \quad \rho(1,x)=\rho_1(x) \, x \in \Omega
\end{array} \right.
\end{equation}
The velocity function $v$, is determined in order to minimize the
functional.
\begin{equation}\label{optiproblem2}
\inf_{\rho,v} \int_0^1 \int_\Omega \rho(t,x) \| v(t,x) \|^2 \, dtdx.
\end{equation}
Thus we get an image sequence through the gray-value map $\rho$.
Let us mention \cite{Tannenbaum}, for example,  where the optimal
mass transportation approach is used in images processing. In this work, 
the optimal transport problem in 2D is decomposed in several 1D optimal transport problems which are
easier to numerically solved. The algorithm proposed here, is
based on a least squares formulation for the transport equation (\cite{Bochev} for example), which
differs from the methods proposed in \cite{Tannenbaum}, or in \cite{Benamou-Brenier}.
Notice that the proposed transport in this paper differs from  the
optimal transportation, studied e.g. in the  book of  C. Villani \cite{Vil07}.
\section{Algorithm for solving the optical flow}\label{algo}
In this section, an algorithm is presented to solve the optical flow given by equations
\eqref{optiproblem}, and \eqref{optiproblem2}.
Let us first specify our hypotheses.
\begin{itemize}
\item [H1] The domain $\Omega$ is a bounded $C^{2,\alpha}$ domain satisfying the
exterior sphere condition.
\item [H2]  The functions $\rho_i \in C^{1,\alpha}(\overline \Omega)$ for $i=0,1$,  with
$\rho_0 = \rho_1 \text{ on }\partial \Omega$.
Moreover there exist two constants such that
$0<\underline \beta \le \rho_i \le \overline \beta$ in $\Omega$.
\end{itemize}
Let $\rho^0 \in C^{1,\alpha}([0,1]\times\overline \Omega)$
be given by $\rho^0(t,x)=(1-t)\rho_0(x) + t \rho_1(x)$. We have
$\|\partial_t \rho^0\|_{C^{0,\alpha}([0,1]\times\overline \Omega)}\le C(\rho_0, \rho_1)$
and $\partial_t \rho^0\vert_{\partial \Omega}=0$.

For each $t\in [0,1]$, our need to solve problem \eqref{optiproblem}-\eqref{optiproblem2} is a
velocity field vanishing on $\partial \Omega$. To do so, the following method is used.

Assume that $\rho^n \in C^{1,\alpha}([0,1]\times\overline \Omega)$ is given with
$\rho^n(0,x)=\rho_0(x)$, and $\rho^n(1,x)=\rho_1(x)$ for all $x \in \Omega$.
\begin{itemize}
\item Compute $\eta^{n+1}$ such that
\begin{equation}\label{problemC}
\left \{ \begin{array}{l}
 -\divg(\rho^{n}(t,\cdot)\nabla \eta^{n+1}) = 0 \, \text{ in }
 \Omega\\
 \displaystyle \rho^{n}(t,\cdot)\derp{\eta^{n+1}}{n}=1 \quad \text{ on }
 \partial \Omega,
\end{array} \right.
\end{equation}
and set
$\displaystyle C^{n+1}(t) = \frac{1}{\vert \partial \Omega \vert} \int_\Omega \partial_t \rho^{n}
\eta^{n+1} \,
dx$.
\item For each $t\in [0,1]$ compute $\varphi^{n+1}$  solution of
\begin{equation}\label{problemphi}
\left \{ \begin{array}{l}
 -\divg(\rho^{n}(t,\cdot)\nabla \varphi^{n+1}) = \partial_t \rho^{n}(t,\cdot), \, \text{ in }
 \Omega\\
\varphi^{n+1}=C^{n+1}(t) \quad \text{ on }
 \partial \Omega.
\end{array} \right.
\end{equation}
\item Set $v^{n+1}= \nabla \varphi^{n+1}$.
\item Compute $\rho^{n+1}$,  {$L^2$-least squares} solution of
\begin{equation}\label{problemrho}
\left \{ \begin{array}{l}
\partial_t \rho^{n+1}(t,x)  + \prodscal{v^{n+1}(t,x)}{\nabla \rho^{n+1}(t,x)} = 0, \;
(t,x) \in (0,1)\times \Omega ,\\
\rho^{n+1}(0,x)=\rho_0(x); \quad \rho^{n+1}(1,x)=\rho_1(x) \; x \in \Omega.
\end{array} \right.
\end{equation}
\end{itemize}
\begin{rem}
\label{vit_norm}
If the requirement  $\rho_{0}{\big \vert_{\partial \Omega}}=\rho_{1}{\big \vert_{\partial \Omega}}$
is canceled in hypothesis H2, then the boundary condition of problem (\ref{problemphi}) is replaced
by $ \displaystyle \derp{\varphi^{n+1}}{n}=0$ and we do not need anymore the constants $C^{n+1}$.
\end{rem}
For each $t \in [0,1]$, since $\rho^n(t,\cdot)$, and $\partial_t \rho^n (t,\cdot)\in
C^{0,\alpha}(\overline \Omega)$, theorem 6.14 p. 107 of \cite{Gilbarg} applies, and
there exists a unique $\varphi^{n+1}(t,\cdot) \in C^{2,\alpha}(\overline \Omega)$ solution of
problem \eqref{problemphi}.
In problem \eqref{problemphi} the  time $t$ is a parameter.
Since $\rho^n \in C^{1,\alpha}$, and $\partial_t
\rho^n \in C^{0,\alpha}$, and $C^{n+1} \in C^{0,\alpha} $, the classical $C^{2,\alpha}(\overline
\Omega)$ a priori estimates for solutions to elliptic problems allow us to show that $\varphi^{n+1}$
is a $C^{0,\alpha}$ function with respect to time.  So we have
$$
\|\varphi^{n+1}\|_{C^{0,\alpha}([0,1];C^{2,\alpha}(\overline \Omega))}\le M(\|C^n
\|_{C^{0,\alpha}([0,1])} + \|\partial_t \rho^n\|_{C^{0,\alpha}([0,1]\times\overline \Omega)}).
$$
Consider the extension of $\varphi^{n+1}$ by $C^{n+1}$ outside of the domain $\Omega$, still denoted
by$\varphi^{n+1}$.
Since the right hand side of equation \eqref{problemphi} vanishes on
$\partial\Omega$, this extension is regular, so the function $v^{n+1}$ vanish
outside $\Omega$ and belongs to $C^{0,\alpha}([0,1];C^{1,\alpha}(\IR^2))$.

Define the flow
$X^{n+1}_{+}(s,t,x)\in C^{1,\alpha}([0,1]\times[0,1]\times\IR^2 ; \IR^2) $ by
\begin{equation}\label{problemflo}
\left \{ \begin{array}{l}
\displaystyle \frac{d}{ds}  X^{n+1}_{+}(s,t,x)  = + v^{n+1}(s,X^{n+1}_{+}(s,t,x))\, \text{ in }
(0,1)\\
X^{n+1}_{+}(t,t,x)  = x.\\
\end{array} \right.
\end{equation}
This flows will be constant for $x \in \partial \Omega$.
Observe that solving the transport equation (\ref{problemrho}) in the least-squares sense is
equivalent to solve this equation on each integral curves defined by (\ref{problemflo}).

Set $r(s)= \rho^{n+1}(s,X^{n+1}_{+}(s,t,x))$,  and  express equation \eqref{problemrho} along
the integral curves of equation \eqref{problemflo}. The equation is reduced to the following
ordinary differential equation with initial and final conditions.
\begin{equation}\label{lsode0}
\left \{ \begin{array}{l}
\displaystyle \frac{d}{ds}  r(s)   = 0\, \\
r(0) = \rho_0(X^{n+1}_{+}(0,t,x)); \, r(1)=\rho_1(X^{n+1}_{+}(1,t,x)).\\
\end{array} \right.
\end{equation}
The $L^2$ least squares solution of \eqref{lsode0} minimizes $\displaystyle \frac{d}{ds}  r(s)$, and is given by:
\begin{equation*}\label{repres0}
r(s)= (1-s) \rho_0(X^{n+1}_{+}(0,t,x)) + s \rho_1(X^{n+1}_{+}(1,t,x)).
\end{equation*}
Therefore the following representation formula for the function $\rho^{n+1}$ is proved.
\begin{lemm}\label{existrho}
The $L^2$-least squares solution of problem \eqref{problemrho} is given by
\begin{equation}\label{repform}
\begin{array}{c}
\rho^{n+1}(t,x)= (1-t) {\rho_0(X^{n+1}_{+}(0,t,x))} +
t {\rho_1(X^{n+1}_{+}(1,t,x))}.
\end{array}
\end{equation}
\end{lemm}
Remark that the regularity of the function $\rho^{n+1}$ is a consequence of the regularity of the
flow $X_+^{n+1}$.

Let us now consider the convergence of the algorithm \eqref{problemC}-\eqref{problemrho}.
\begin{theo}\label{convalgo}
There exist
$(\rho,\varphi)\in C^{1}([0,1]\times\overline \Omega; \IR_+^*)\times
C^{0}([0,1];C^{2}(\overline \Omega))$,
$L^2$-least squares solution, respectively solution of
\begin{equation}\label{problemlim1a}
\left \{ \begin{array}{l}
\partial_t \rho  + \prodscal{\nabla \varphi}{\nabla \rho} = 0, \quad \mathrm{in} \,
(0,1)\times \Omega\\
\rho(0,x)=\rho_0(x); \quad \rho(1,x)=\rho_1(x) \quad \mathrm{in} \, \Omega \\
\end{array} \right.
\end{equation}
\begin{equation}\label{problemlim1b}
\left \{ \begin{array}{l}
 -\divg(\rho(t,\cdot)\nabla \varphi) = \partial_t \rho(t,\cdot), \quad \mathrm{in} \,
 \Omega\\
\varphi=C(t) \quad \mathrm{and}  \quad \nabla  \varphi = 0 \quad \mathrm{on} \,
 \partial \Omega\\
\end{array} \right.
\end{equation}
with $C(t)$ defined as follow.
\begin{equation}\label{problemlim2}
\left \{ \begin{array}{l}
 -\divg(\rho(t,\cdot)\nabla \eta) = 0 \; \mathrm{in} \, \Omega\\
\rho(t,\cdot) \, \partial_n \eta=1 \; \mathrm{on} \,
 \partial \Omega\\
 C = \displaystyle \frac{1}{\vert \partial \Omega \vert} \int_\Omega \partial_t \rho \, \eta \, dx.
\end{array} \right.
\end{equation}
\end{theo}
\begin{dem}
Since
$\|v^0 \|_{C^{0,\alpha}([0,1])} + \|\partial_t \rho^0\|_{C^{0,\alpha}([0,1]\times\overline
\Omega)}$
is bounded,
$\|\varphi^{n+1}\|_{C^{0,\alpha}([0,1];C^{2,\alpha}(\overline \Omega))}$
and $\|v^{n+1}\|_{C^{0,\alpha}([0,1];C^{1,\alpha}(\IR^2))}$ are uniformly bounded in $n$.

From lemma \ref{existrho} there exists a unique
$\rho^{n+1}$, the $L^2$-least squares solution of \eqref{problemrho}.
Let us give an estimate for $D_3X_{+}^{n+1}$.
Starting from
$$
D_1 X^{n+1}_{+}(s,t,x))= v^{n+1}(s,X^{n+1}_{+}(s,t,x)),
$$
we deduce (see \cite{Ambrosio})
\begin{equation}
\left \{ \begin{array}{l}
D_3D_1 X^{n+1}_{+}(s,t,x)  = D_2 v^{n+1}(s,X^{n+1}_{+}(s,t,x))D_3X^{n+1}_{+}(s,t,x)\\
D_3X^{n+1}_{+}(t,t,x)  = Id.\\
\end{array} \right.
\end{equation}
Since $D_3D_1X^{n+1}_{+}(s,t,x)  =D_1D_3X^{n+1}_{+}(s,t,x)$ we get
\begin{equation}\label{repDX}
D_3X^{n+1}_{+}(s,t,x)  = e^{-\int_t^s D_2(v^{n+1}(\tau,X^{n+1}_{+}(\tau,t,x)))\, d\tau} Id.
\end{equation}
Thus $\|D_3v_{+}^{n+1}\|_{C^{0,\alpha}([0,1]^2\times \IR^2)}$ is uniformly bounded in $n$.
Moreover we have \cite{Ambrosio}
$$
D_2X^{n+1}_{+}(s,t,x)= \prodscal{v^{n+1}(s,t,x)}{D_3X^{n+1}_{+}(s,t,x)}
$$
so $\|D_2v^{n+1}\|_{C^{0,\alpha}([0,1]^2\times \IR2)}$ is bounded independently of $n$.

From theorem \ref{existrho} we deduce that
$\|\rho^{n+1}\|_{C^{1,\alpha}([0,1]\times  \overline \Omega)}$ is uniformly bounded. Since the
embeddings
$$
C^{0,\alpha}([0,1];C^{2,\alpha}(\overline \Omega))\hookrightarrow
C^{0}([0,1];C^{2}(\overline \Omega))\, \mathrm{and} \,  C^{1,\alpha}([0,1]\times  \overline
\Omega)\hookrightarrow C^{1}([0,1]\times  \overline \Omega)
$$
are relatively compact there is a subsequence of $(\rho^{n},\varphi^{n})$
solution of \eqref{problemC}-\eqref{problemrho}, still denoted by $(\rho^{n},\varphi^{n})$
converging to
$(\rho,\varphi)$ in $C^{1}([0,1]\times  \overline \Omega)\times C^{0}([0,1];C^{2}(\overline
\Omega))$, and $(\rho,\varphi)$ is the solution of
\eqref{problemlim1a}-\eqref{problemlim2} provided the boundary conditions to be  justified.
The condition $\nabla\varphi^n\vert_{\partial \Omega} =0$ is valid for the approximations
$\varphi^n$ (since the functions can be extended by $C^n$ outside of $\Omega$). So the convergence
in $C^{0}([0,1];C^{2}(\overline \Omega))$ yields the condition for the gradient of limit function.
For the approximations of function $\rho$, the formula given in lemma \ref{existrho} combined with
the regularity result show that the boundary conditions are exactly satisfied. These conditions are
thus valid for the limit function due to the convergence in $C^1$.
\end{dem}
We will show in the next section that the solution minimizes the same energy as for
the time dependent optimal transportation mass.
\section{Interpretation of solutions to problem \eqref{problemlim1a}-\eqref{problemlim2}}
\label{masse}
\noindent In this section it is shown that the solution to problem
\eqref{problemlim1a}-\eqref{problemlim2} is a solution to a time dependent
 mass transportation problem. 

Zero is a bound from below of the following functional $\frac{1}{4}\int_0^1 \| \partial_t u  +
\divg(u \nabla (\psi-C) ) \|_{L^{2}(\Omega)}^2 \, dt $ to be minimized with respect to
$(u,\psi)$.
$$ 0= \frac{1}{4}\int_0^1 \| \partial_t \rho  +
\divg(\rho \nabla (\varphi-C) ) \|_{L^{2}(\Omega)}^2 \, dt
$$
thus $(\rho,\varphi-C)$ solution of \eqref{problemlim1a}-\eqref{problemlim2} minimizes
\begin{equation*}
\underset{
\begin{array}{c}
\scriptstyle \{\psi\in L^2((0,1);H^1_0(\Omega)), \; u\in L^2((0,1);L^2(\Omega)) \\
\scriptstyle \partial_t u  + \prodscal{u} {\nabla \psi} = 0 \\
\scriptstyle u(0)=\rho_0;\; u(1)=\rho_1 \text{ in } \Omega \}
\end{array}
}
{\rm Min}
\frac{1}{4} \int_0^1 \| \partial_t u  +
\divg(u \nabla \psi ) \|_{L^2(\Omega)}^2 \, dt.
\end{equation*}
So the solution $(\rho,\varphi)$ of problem
\eqref{problemlim1a}-\eqref{problemlim2}, satisfies
\begin{equation}\label{problemint1}
(\rho,\varphi-C)=
\underset{
\begin{array}{c}
\scriptstyle \{\psi\in L^2((0,1);H^1_0(\Omega)), \\
\scriptstyle u\in L^2((0,1);L^2(\Omega)), \\
\scriptstyle \partial_t u  + \prodscal{u} {\nabla \psi}  = 0, \\
\scriptstyle u(0)=\rho_0;\; u(1)=\rho_1 \text{ in } \Omega \}
\end{array}
}
{\text {Argmin} }
\frac{1}{4} \int_0^1 \| \partial_t u  +
\divg(u \nabla \psi ) \|_{H^{-1}} \, dt.
\end{equation}
\begin{theo}\label{problemwas}
Let $(\rho,\varphi)$ be the solution of problem
 \eqref{problemlim1a}-\eqref{problemlim2} given in
 theorem \ref{convalgo}, then it satisfies
\begin{equation}\label{eqwas}
(\rho,\nabla\varphi)=
\underset{
\begin{array}{c}
\scriptstyle \{v\in L^2((0,1);\left (H^1(\Omega)\right )^2), \\
\scriptstyle u\in L^2((0,1);L^2(\Omega)), \\
\scriptstyle \partial_t u  + \divg(uv)=0, \\
\scriptstyle \partial_t u  + \prodscal{v}{\nabla u} =0, \\
\scriptstyle u(0)=\rho_0; \, u(1)=\rho_1
\text{ in } \Omega \}
\end{array}
}
{\rm Argmin}
\int_0^1 \int_\Omega u \| v \|^2 \, dxdt.\\
\end{equation}
\end{theo}
\begin{dem} Observe that in problem \eqref{eqwas}, the transport equation is solved with a
$L^2$-least square procedure. This is equivalent to find $\xi= u -(1-t)\rho_0 + t \rho_1$ such that
$$
\partial_t \xi  + \prodscal{v}{\nabla \xi} = P_R( \rho_1-\rho_0 +
\prodscal{v}{\nabla \left ((1-t)\rho_0 + t \rho_1\right )}
$$
where $P_R$ is the $L^2$ projection onto the range of the transport operator for the velocity $v$.

For $u \in L^\infty(\Omega)$, the expression  $v \mapsto \divg(uv)$ is well
defined as a linear continuous on  $H^1_0(\Omega)$. If moreover $0<\underline \beta \le u \le
\overline \beta$ in $\Omega$,
let $H=H^1_0(\Omega)$ be equipped with the following inner product.
$$
(\theta,\psi)=\int_\Omega u \prodscal{\nabla \theta}{\nabla \psi} \, dx,
$$
which induces a semi-norm which is equivalent to the $H^1$-norm.
The Riesz theorem claims that for the linear continuous form
$$
\mathcal{L}_v(\psi) =<-\divg{(uv)},\psi>_{H;H'},
$$
there is a unique $\theta \in H$ such that
$$
\mathcal{L}_v(\psi)=\int_\Omega u \prodscal{\nabla \theta}{\nabla \psi} \, dx,
\, \forall \psi \in H.
$$
Therefore $v=\nabla \theta$ for a $\theta \in H^1_0(\Omega)$ and problem \eqref{eqwas} is
reduced to
\begin{equation}\label{probleminfv2}
(\rho,\nabla\varphi)=
\underset{
\begin{array}{c}
\scriptstyle \{ \psi \in L^2((0,1);H^1_0(\Omega)\cap H^2(\Omega) ), \\
\scriptstyle u\in L^2((0,1);L^2(\Omega)), \\
\scriptstyle \partial_t u  + \divg(u \nabla \psi)=0, \\
\scriptstyle \partial_t u  +\prodscal{\nabla \psi}{u}=0,\\
\scriptstyle \, u(0)=\rho_0; \, u(1)=\rho_1 \text{ in } \Omega \}
\end{array}
}
{\rm Argmin}
  \int_0^1 \int_\Omega u \| \nabla \psi \|^2 \, dxdt.
\end{equation}
or
\begin{equation}\label{probleminfv2bis}
(\rho,\nabla\varphi)=
\underset{
\begin{array}{c}
\scriptstyle \{-\psi \in L^2((0,1);H^1_0(\Omega)\cap H^2(\Omega)), \\
\scriptstyle u\in L^2((0,1);L^2(\Omega)), \\
\scriptstyle \partial_t u  - \divg(u \nabla \psi)=0, \\
\scriptstyle \partial_t u  - \prodscal{u} {\nabla \psi}=0, \\
\scriptstyle u(0)=\rho_0; \, u(1)=\rho_1 \text{ in } \Omega \}
\end{array}
}
{\rm Argmin}
  \int_0^1 \int_\Omega u \| \nabla \psi \|^2 \, dxdt.
\end{equation}
Since
$$
\int_\Omega u \| \nabla \psi \|^2 \, dx = \|-\divg(u\nabla \psi)\|^2_{H^{-1}},
$$
problem \eqref{probleminfv2} reads:
\begin{equation}\label{probleminfv3}
(\rho,\nabla\varphi)=
\underset{
\begin{array}{c}
\scriptstyle \{-\psi \in L^2((0,1);H^1_0(\Omega)\cap H^2(\Omega)), \\
\scriptstyle u\in L^2((0,1);L^2(\Omega)), \\
\scriptstyle \partial_t u  - \divg(u \nabla \psi)=0, \\
\scriptstyle \partial_t u  - \prodscal{u} {\nabla \psi}=0, \\
\scriptstyle u(0)=\rho_0; \, u(1)=\rho_1 \text{ in } \Omega \}
\end{array}
}
{\rm Argmin}  \int_0^1 \|\divg(-u\nabla \psi)\|^2_{H^{-1}} \, dt.
\end{equation}
From the definition of the linear form $\mathcal L$, observe that
$$
\frac{1}{4} \|\divg(-u\nabla \psi) + \partial_t u\|^2_{H^{-1}} =
\|\divg(-u\nabla \psi)\|^2_{H^{-1}},
$$
so problem \eqref{probleminfv3} reads:
\begin{equation}\label{probleminfv4}
(\rho,\nabla\varphi)=
\underset{
\begin{array}{c}
\scriptstyle \{\psi \in L^2((0,1);H^1_0(\Omega)\cap H^2(\Omega)), \\
\scriptstyle u\in L^2((0,1);L^2(\Omega)), \\
\scriptstyle \partial_t u  + \divg(u \nabla \psi)=0, \\
\scriptstyle \partial_t u  + \prodscal{u} {\nabla \psi}=0, \\
\scriptstyle u(0)=\rho_0; \, u(1)=\rho_1 \text{ in } \Omega \}
\end{array}
}
{\rm Argmin}  \frac{1}{4}\int_0^1 \|\divg(u\nabla \psi) + \partial_t u\|^2_{H^{-1}}\, dt.
\end{equation}
Gathering \eqref{problemint1} with the previous result proves the theorem.
\end{dem}
\begin{rem}
The proposed transport is a one which minimizes the divergence of the velocity in a
weighted dual norm of $H^1_0$
\end{rem}
\begin{rem}
The Dacorogna-Moser transport \cite{dacorogna}  $\rho(t,x)=(1-t)\rho_0(x) + t \rho_1(x)$ with
\begin{equation}\left \{ \begin{array}{l}
- \Delta \varphi = \partial_t \rho \text{ in } \Omega; \\
\partial_n \varphi = 0; \text{ on } \partial \Omega;
\end{array}\right.
\end{equation}
and the velocity $v =\frac{\nabla \varphi}{\rho}$ satisfies
\begin{equation}\label{eqDaco}
(\rho,\nabla\varphi)=
\underset{
\begin{array}{c}
\scriptstyle \{v\in L^2((0,1);\left (H^1(\Omega)\right )^2), \\
\scriptstyle u\in L^2((0,1);L^2(\Omega)), \\
\scriptstyle \partial_t u  + \divg(v) =0,\\
\scriptstyle u(0)=\rho_0; \, u(1)=\rho_1 \text{ in } \Omega\}
\end{array}
}
{\rm Argmin}
\int_0^1 \int_\Omega  \| v \|^2 \, dxdt.\\
\end{equation}
\end{rem}
Indeed, as before, consider the space $H=H^1_0$ equipped with the previous
semi-norm, and set
$$
\mathcal{L}_v(\psi) =<-\divg{(v)},\psi>_{H;H'}=<\partial_tu,\psi>_{H;H'}.
$$
Then the problem is reduced to
\begin{equation}\label{probleminfv6}
(\rho,\nabla\varphi)=
\underset{
\begin{array}{c}
\scriptstyle \{\psi \in L^2((0,1);H^1_0(\Omega)\cap H^2(\Omega)), \\
\scriptstyle u\in L^2((0,1);L^2(\Omega)) \\
\scriptstyle\partial_t u  + \divg( \nabla \psi )=0 \\
\end{array}
}
{\rm Argmin} \frac{1}{4} \int_0^1 \|\partial_t u - \divg(\nabla \psi)\|^2_{H^{-1}} \, dt.
\end{equation}
Using the relation $\partial_t u = -\divg(\nabla \psi)$ and Jensen's inequality we get
$$
 \int_0^1 \|\divg(\nabla \psi)\|^2_{H^{-1}} \, dt =
\int_0^1 \|\partial_t u\|^2_{H^{-1}} \, dt \ge \|\int_0^1 \partial_t u  \, dt\|^2_{H^{-1}}= \|
\rho_1 - \rho_0 \|^2_{H^{-1}}.
$$
Thus the Dacorogna-Moser transport is a minimum of the functional \eqref{eqDaco},
so it is a transport minimizing the
divergence of the velocity in $H^{-1}$-norm.
\section{Numerical Approximation of the 2D Optimal Extended Optical Flow}
\label{exemples}
\noindent The numerical method is based on a finite element time-space
$L^2$ least squares formulation (see \cite{Besson}) of the transport problem
\eqref{problemrho}. Define $\vtld^{n+1}$ as
\[\vtld^{n+1} = (1,v^{n+1}_{1},v^{n+1}_{2})^t \]
and for  a sufficiently regular function $\varphi$ defined on $Q$,
set
\[\nabtld\varphi =
\left( \derp{\varphi}{t}, \derp{\varphi}{x_{1}},
\derp{\varphi}{x_{2}} \right)^t,
\]
and
\[
\divtld(\vtld^{n+1}  \ \varphi) = \derp{\varphi}{t} +  \sum_{i=1}^{2}
\derp{ }{x_{i}}( v^{n+1}_{i} \ \varphi).
\]
Let $\{\varphi_1 \cdot \cdot \cdot \varphi_N\}$ be a basis of a space-time finite element
subspace
\[
V_h = \{\varphi, \text{ piecewise regular polynomial functions, with } \varphi(0,\cdot)=
\varphi(1,\cdot)=0 \},
\]
for example, a  brick Lagrange finite element of order one (\cite{Besson2}).
Let $\Pi_h$ be the Lagrange interpolation operator. Let also $W_h$ be the
finite element subspace of $H^1_0(\Omega)$, where the basis functions
$\{\psi_1 \cdot \cdot \cdot \psi_M\}$ are  the traces at $t=0$ of basis functions
$\{\varphi_i\}_{i=1}^N$.
An approximation of problem \eqref{problemphi} is the following. For a discrete sequence of time $t$
compute
\begin{multline}
\int_\Omega(\rho_h^{n}(t,\cdot)\prodscal{\nabla (\varphi_h^{n+1}-C^n(t))}{\nabla
\psi_h}\, dx = \\
\int_\Omega \partial_t \rho_h^{n}(t,\cdot) \psi_h \, dx \quad \forall  \psi_h \in W_h,
\end{multline}
and define  $v^{n+1}= \nabla \varphi_h^{n+1}$.
The $L^2$ least squares formulation of problem \eqref{problemrho} is defined in the following
way. Consider the functional
\[
J(c_h) = \frac{1}{2}  \int_{Q} \prodscal{\vtld^{n+1} }
{\nabtld c_h + \nabtld\left[\Pi_h \big ((1-t)\rho_{0} +
t\rho_{1}\big )\right]}^2 \, dx \, dt.
\]
This functional is convex and coercive in the appropriate anisotropic
Sobolev's space
$H = \{ \varphi \in L^2(Q); \, \prodscal{\vtld^{n+1} }{\nabtld{\varphi}} \in  L^2(Q);
\varphi(0,\cdot)=\varphi(1,\cdot) =0\}$ since the velocity field $v^{n+1} $ is regular enough.
Moreover, $\norml{\prodscal{\vtld^{n+1} }{\nabtld{\varphi}}}{Q} $ is a norm in $H$ (see \cite{Besson}). Set
$$
D=\left [\underline \beta, \overline \beta \right ], \quad \theta_h = \Pi_h \big ((1-t)\rho_{0} +
t\rho_{1}\big ),
$$
and introduce
$$
K_h=\{ \varphi_h \in V_h; \,  \varphi_h + \theta_h \in D \}
$$
it is a closed convex subspace. Thus an approximation of \eqref{problemrho} subject to the
constraint $\rho^{n+1}\in D$ is defined with the following minimization problem:
\begin{equation}\label{minpbcont}
\min_{c \in K_h} J(c).
\end{equation}
The minimizer  of problem \eqref{minpbcont} is
$$
c_h=\rho^{n+1}_h-\theta_h
$$
which is characterized by the Fermat's rule \cite{Aubin}
\begin{multline}\label{dirihomoh}
\int_{Q} \prodscal{\vtld^{n+1} }{\nabtld c_h} \prodscal{\vtld^{n+1} }{\nabtld \psi_h} \, dx \
dt \ge  \\
\int_{Q} \prodscal{-\vtld^{n+1} }{\nabtld \theta_h} \prodscal{\vtld^{n+1} }{\nabtld \psi_h} \, dx \ dt
\end{multline}
for all $\psi_h  \in T_{K_h}$, the contingent cone to $K_h$.
Thus an approximation of the solution to
problem \eqref{problemrho} is
$$
\rho^{n+1}_h =c_h + \Pi_h \big ((1-t)\rho_{0} +
t\rho_{1}\big )\in V_h.
$$
\begin{rem}
Define the bilinear form 
$a_v^{n+1}(\cdot,\cdot)$ on $V_h \times V_h$ by
$$
a_v^{n+1}(\varphi,\psi)= \int_{Q} \prodscal{\vtld^{n+1} }{\nabtld \psi} 
\prodscal{\vtld^{n+1} }{\nabtld \psi_h} \, dx \ dt.
$$
Then the function $c_h$ is the solution of
$$
a_v^{n+1}(c_h,\psi)= a_v^{n+1}(-\theta_h,\psi)
$$
as long as $c_h$ is everywhere non negative.
Moreover, if $V^+_h$  consists of the non negative functions in $V_h$, the function $c_h$ is 
the orthogonal projection of $\theta_h$ onto $V^+_h$, according to the inner product induced by the
bilinear form $a_v^{n+1}$.

When the approximation with finite element of the algorithm proposed in section 
\ref{algo} has converged, the computed solution $\rho_h$ can be decomposed as $\rho_h = \theta_h +
P^{a_v}_{V_h^+} \theta_h$, that is to say the approximated Dacoragna-Moser transport solution
augmented with its orthogonal projection  onto non negative functions space according to the inner
product $a_v$.
\end{rem}

In the next subsections, all the computations are done according to remark \ref{vit_norm}. So
it is assumed that the normal velocity satisfy $\displaystyle v^{n+1}_{n} =
\derp{\varphi^{n+1}}{n}=0$.
\subsection{The transport of a bump}
\label{sec:bosse}
As a first example, the displacement of a bump is considered. More precisely, let $\Omega =
]-1,1[ \times ]-1,1[$, $y_0 = 0.5$, $r_0 = 0.3$, $\beta > 0$, and $\alpha > 0$.
For $(x,y) \in \Omega$, set $r = \sqrt{x^2 +(y-y_0)^2}$, and define
\begin{equation}
\rho_0 = \left\lbrace
\begin{array}{l}
\beta + \alpha e^{-r^2/(r_0^2-r^2)} \quad \text{if} \; r_0^2 > r^2 \\
\beta \quad \text{elsewhere.}
\end{array}
\right.
\end{equation}
Then for $r = \sqrt{x^2 +(y+y_0)^2}$, define
\begin{equation}
\rho_1 = \left\lbrace
\begin{array}{l}
\beta + \alpha e^{-r^2/(r_0^2-r^2)} \quad \text{if} \; r_0^2 > r^2 \\
\beta \quad \text{elsewhere,}
\end{array}
\right.
\end{equation}
The domain $\Omega$ is subdivided into $40 \times 40$ elements, and the time interval $]0,T[ =
]0,1[$ is subdivided into $60$ elements, thus the linear system has $102541$ unknowns. Is is solved
with a preconditioned conjugate gradient.

\begin{figure}[H]
\begin{center}
\begin{minipage}{0.48\linewidth}
\epsfig{file=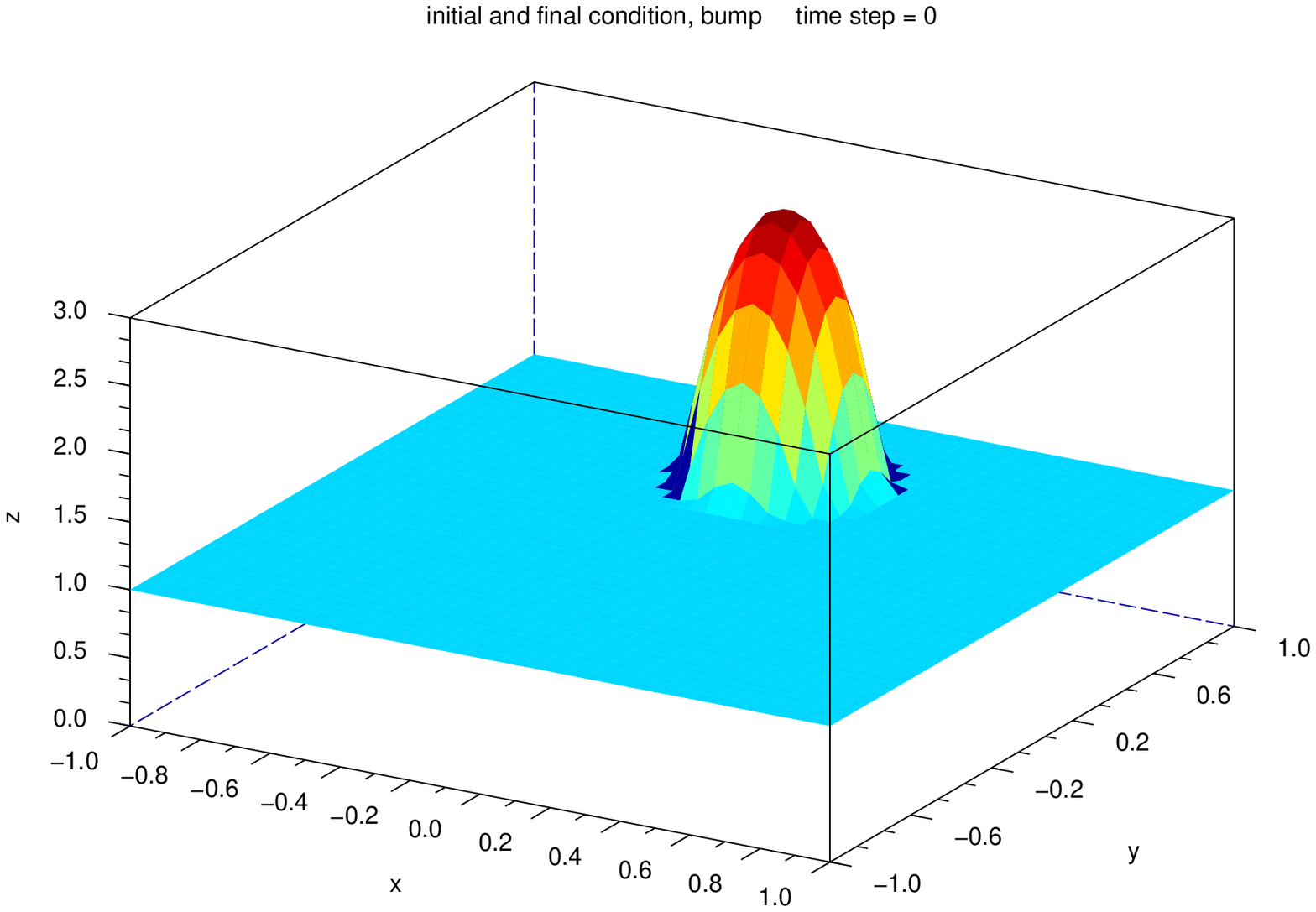,width=7cm}
\end{minipage}
\begin{minipage}{0.48\linewidth}
\epsfig{file=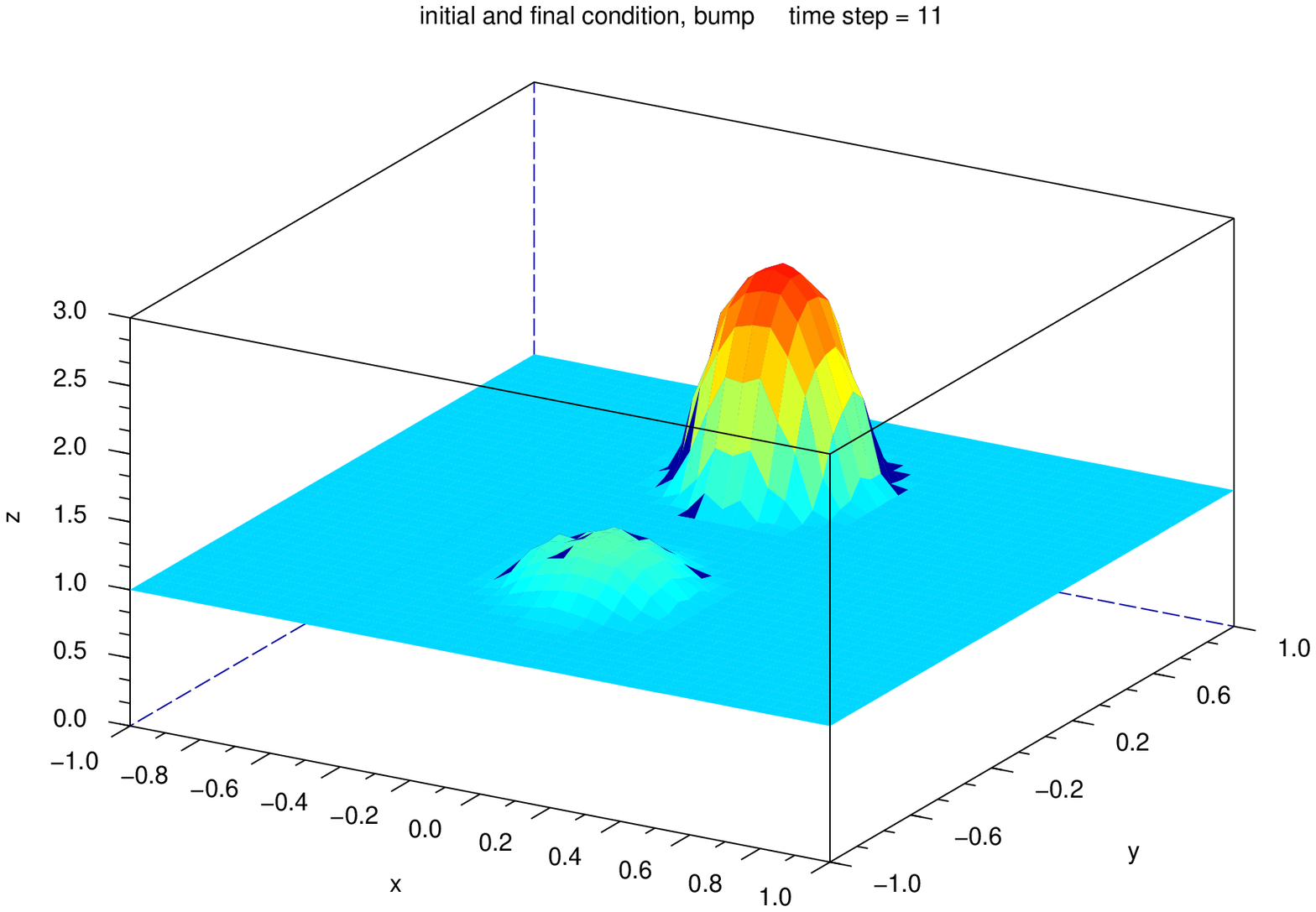,width=7cm}
\end{minipage}
\begin{minipage}{0.48\linewidth}
\epsfig{file=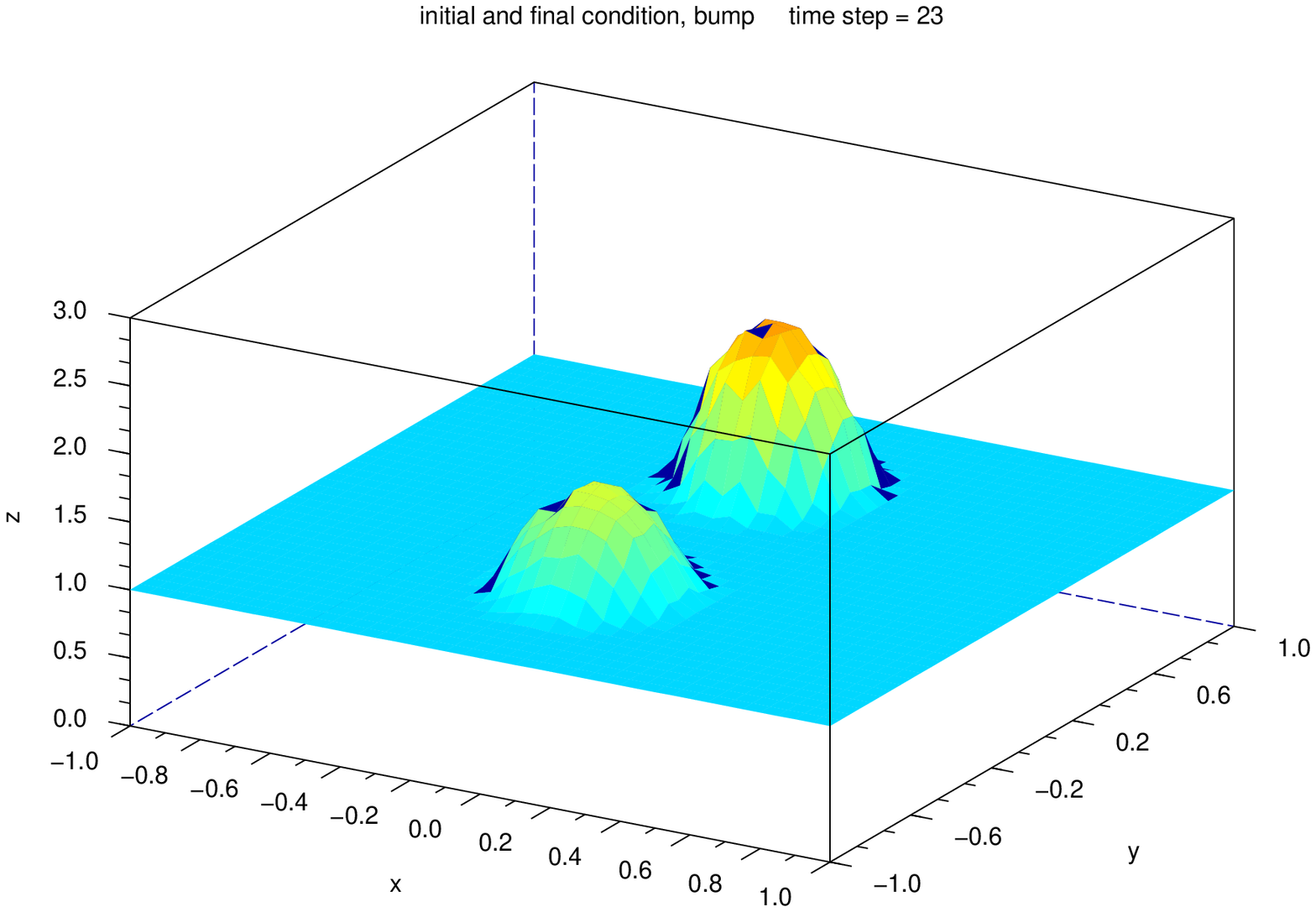,width=7cm}
\end{minipage}
\begin{minipage}{0.48\linewidth}
\epsfig{file=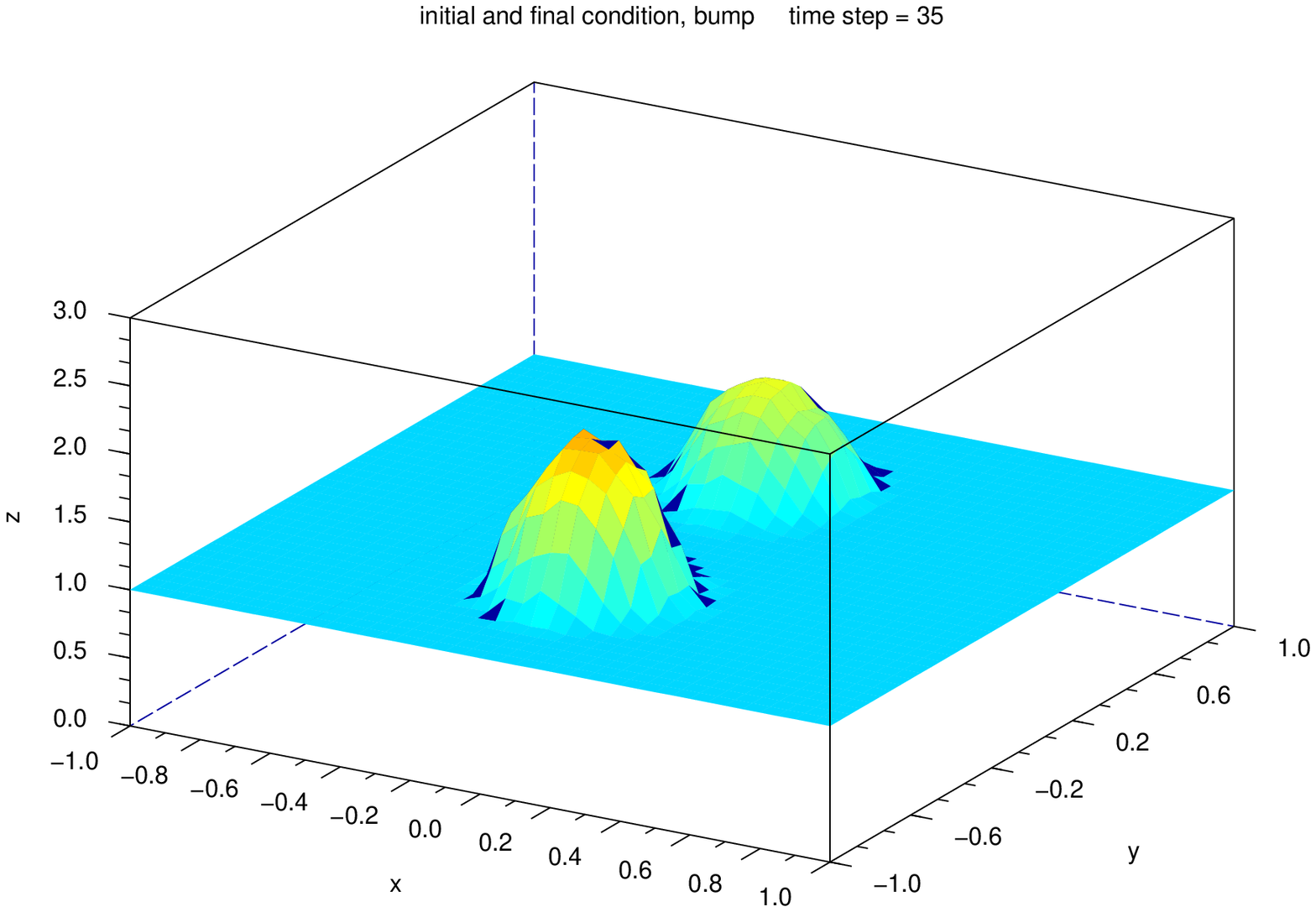,width=7cm}
\end{minipage}
\begin{minipage}{0.48\linewidth}
\epsfig{file=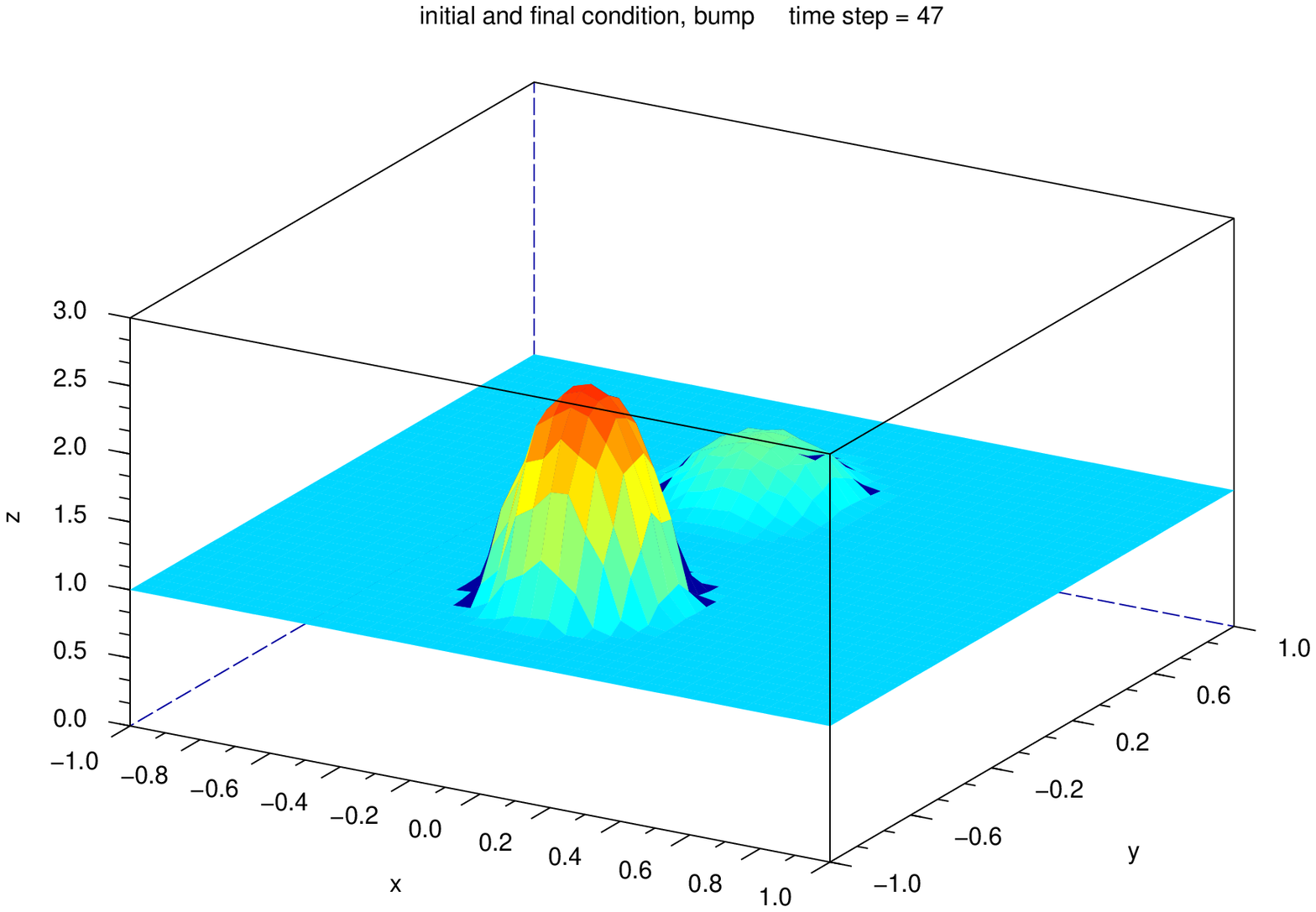,width=7cm}
\end{minipage}
\begin{minipage}{0.48\linewidth}
\epsfig{file=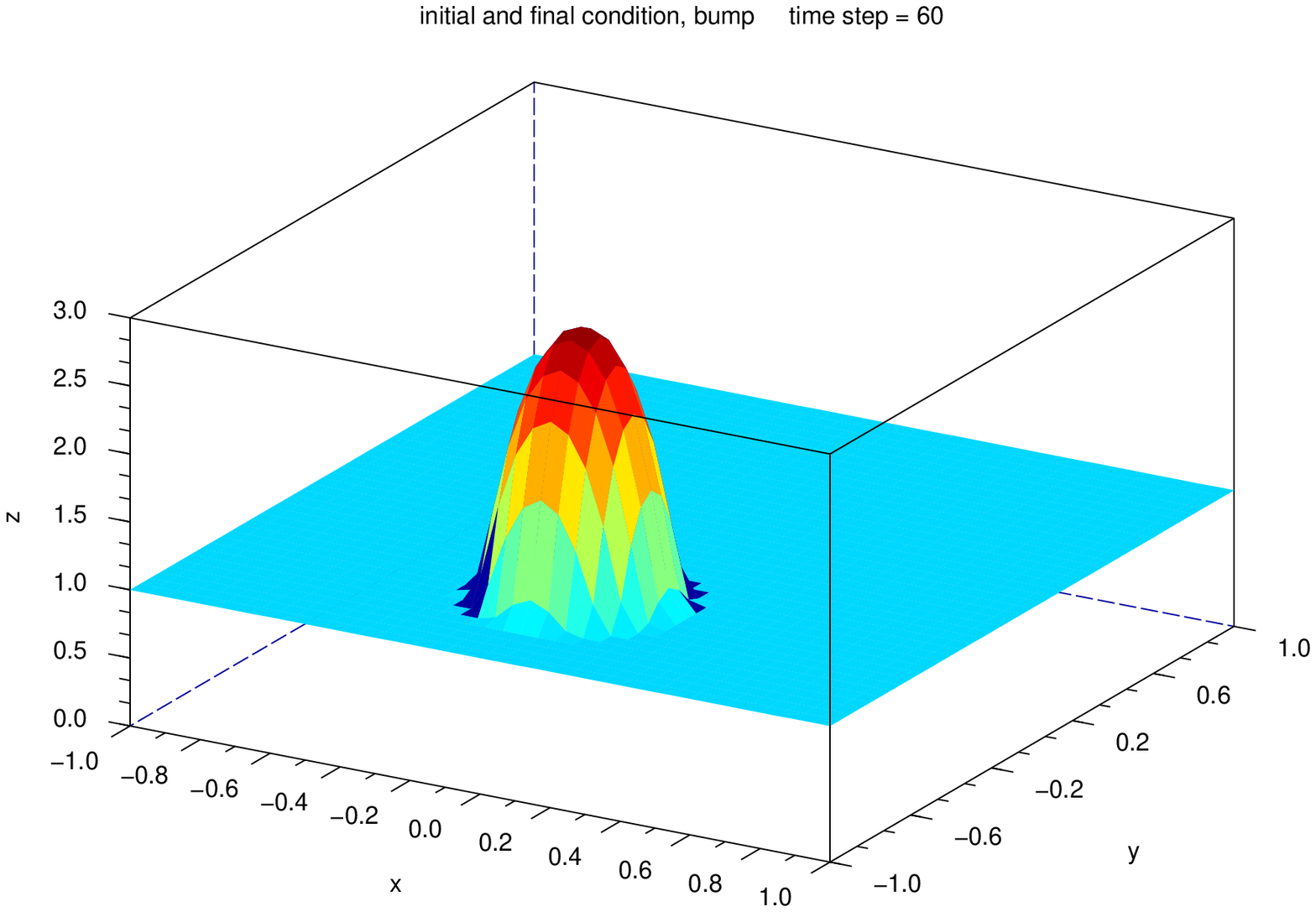,width=7cm}
\end{minipage}
\caption{Bump transport for $\beta=1$}
\label{beta1}
\end{center}
\end{figure}

\begin{figure}[H]
\begin{center}
\begin{minipage}{0.48\linewidth}
\epsfig{file=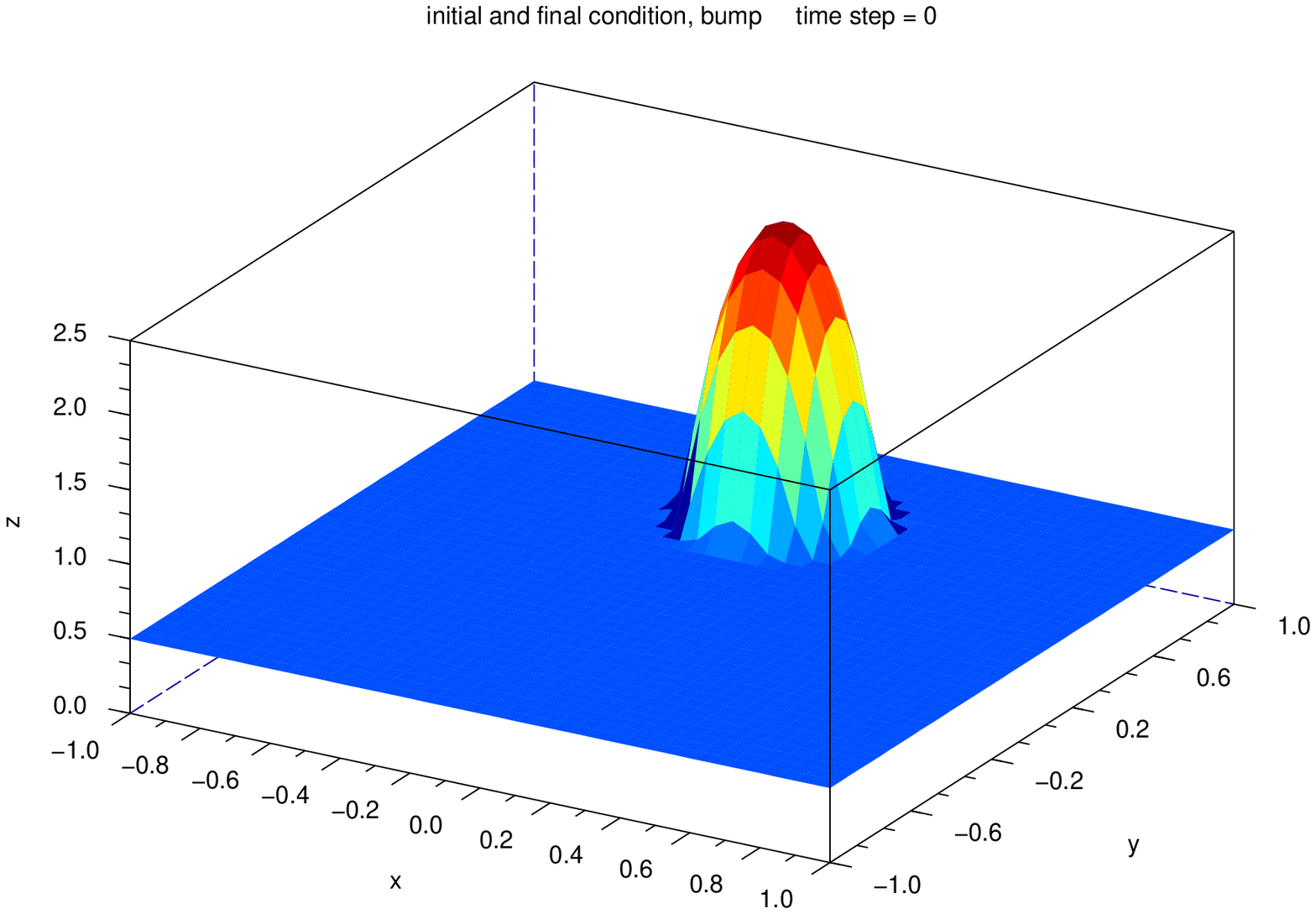,width=7cm}
\end{minipage}
\begin{minipage}{0.48\linewidth}
\epsfig{file=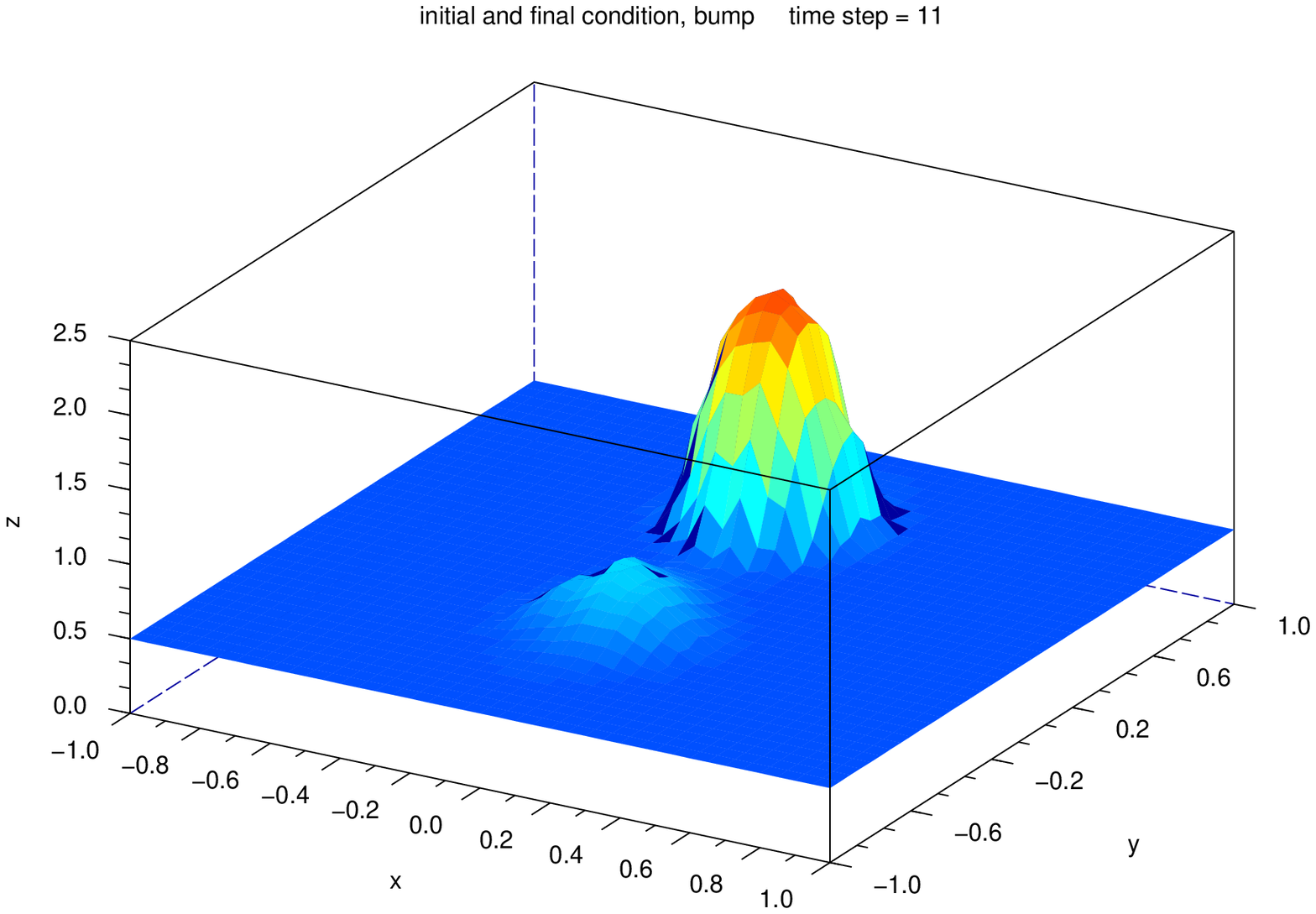,width=7cm}
\end{minipage}
\begin{minipage}{0.48\linewidth}
\epsfig{file=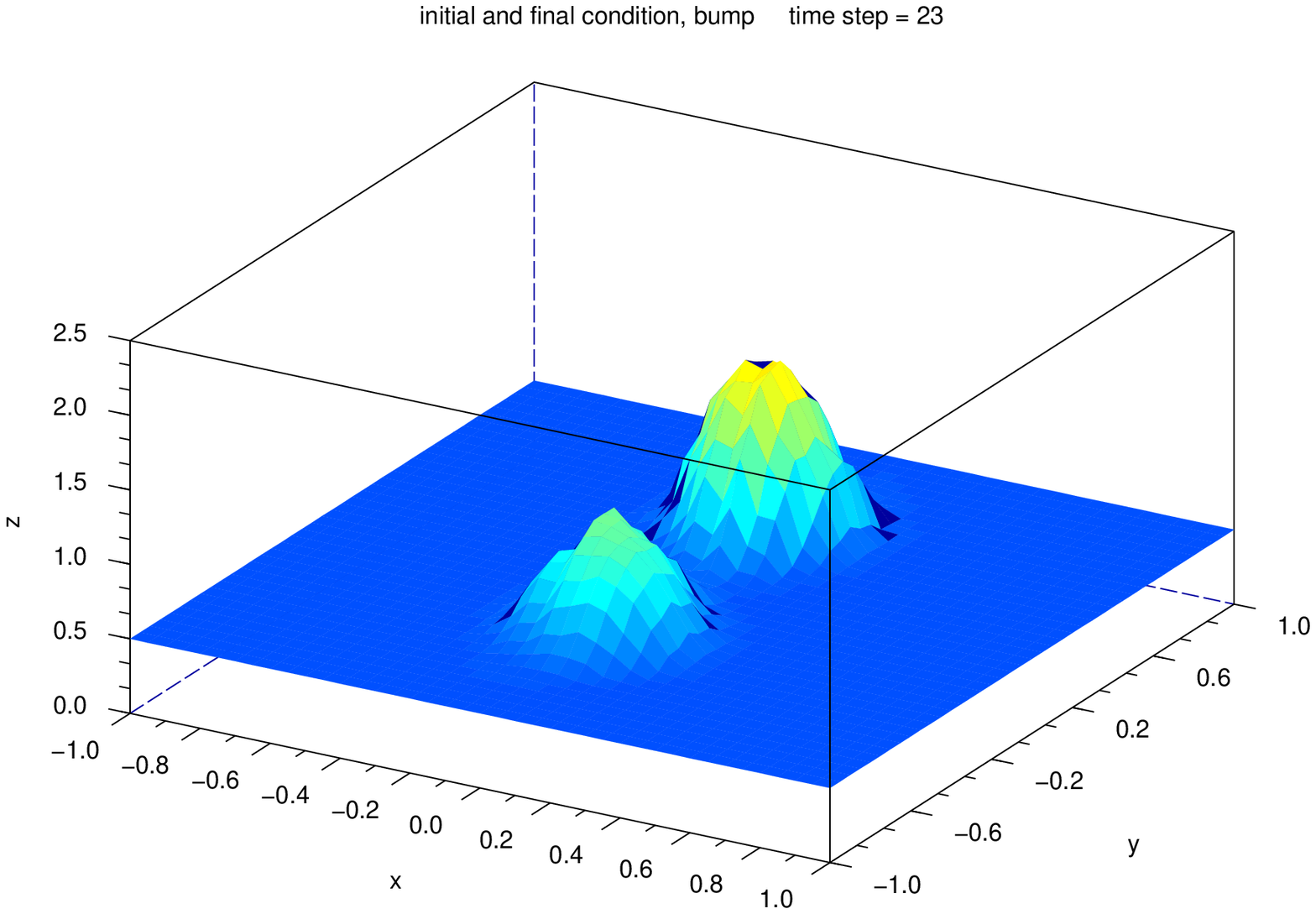,width=7cm}
\end{minipage}
\begin{minipage}{0.48\linewidth}
\epsfig{file=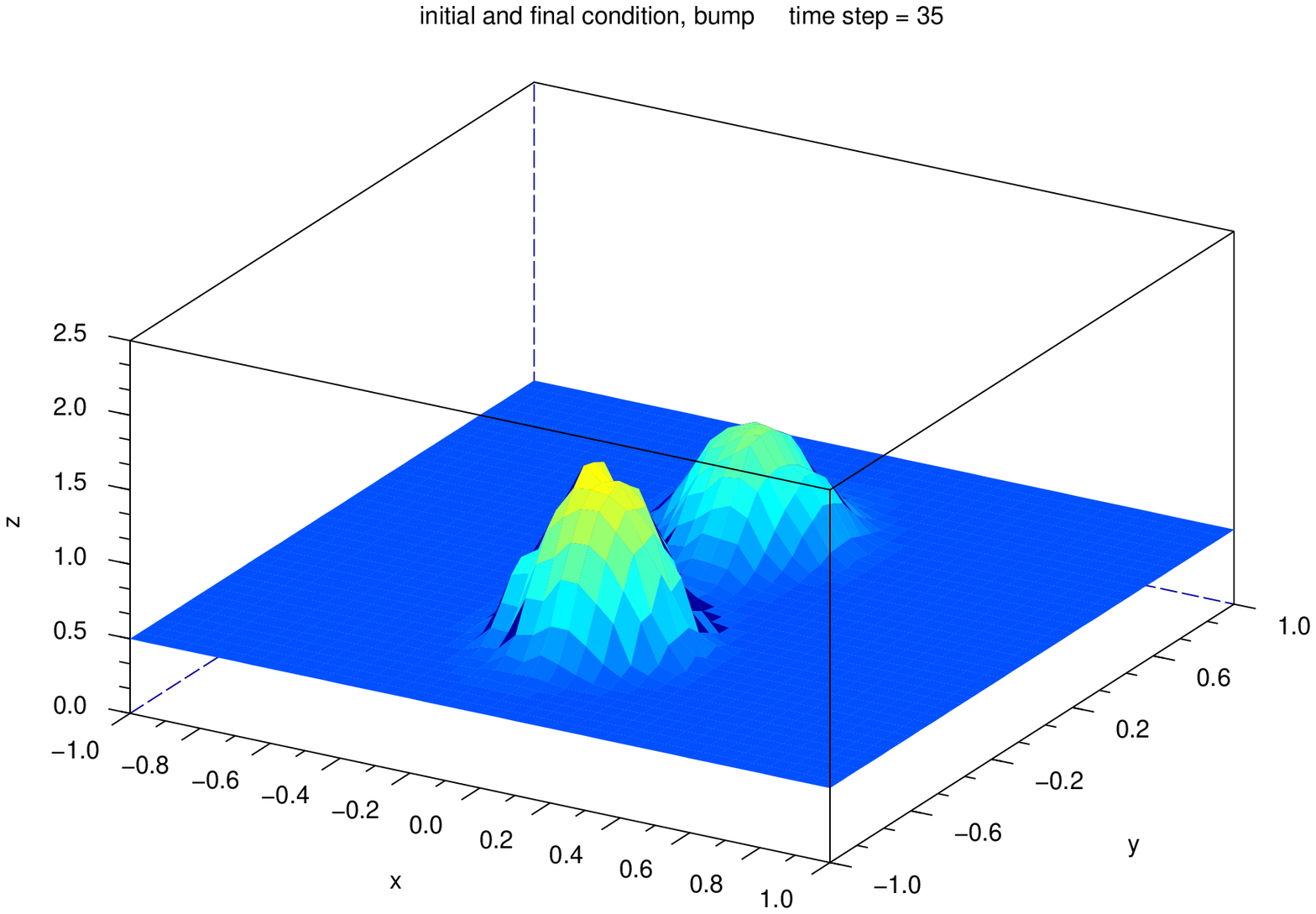,width=7cm}
\end{minipage}
\begin{minipage}{0.48\linewidth}
\epsfig{file=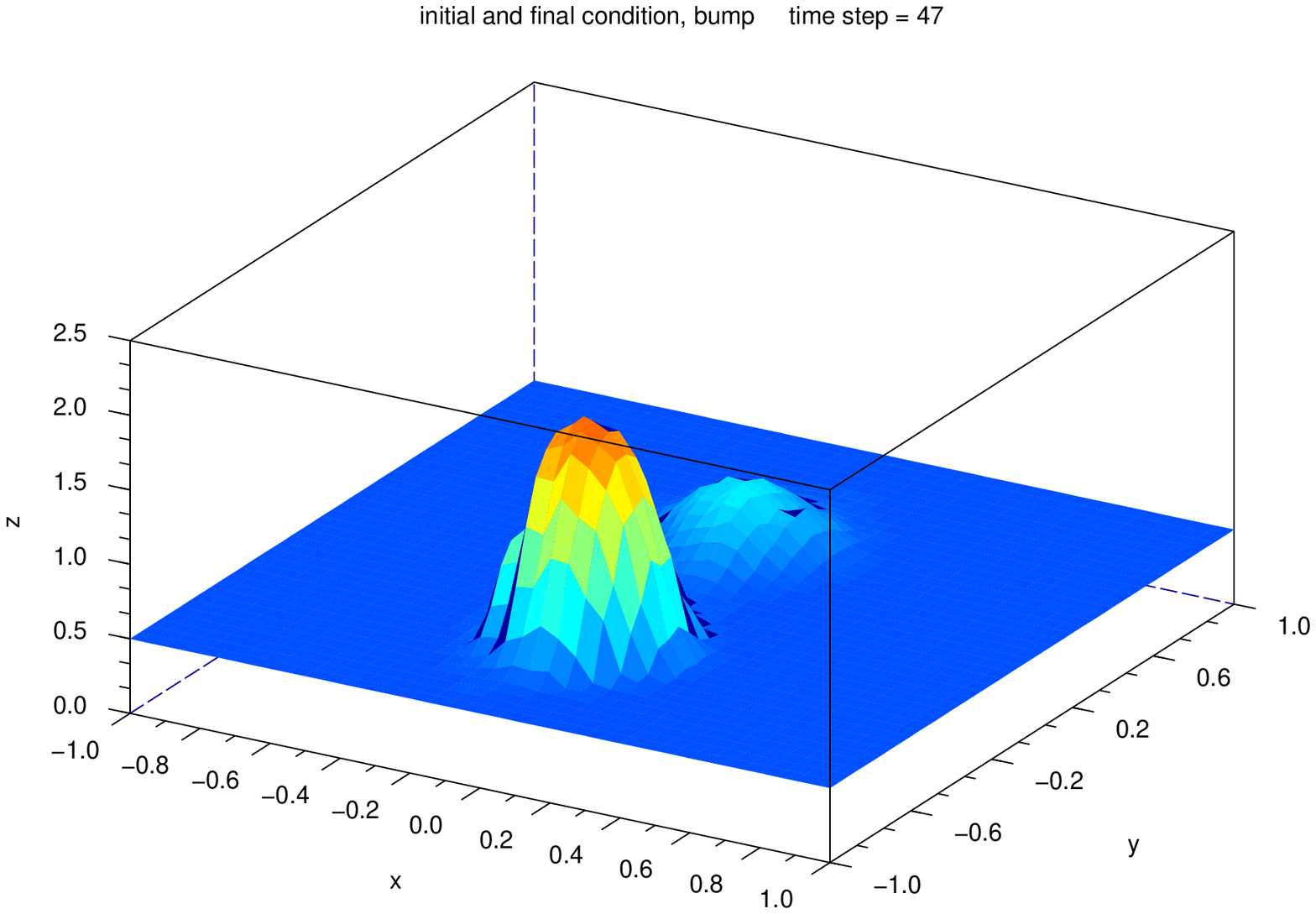,width=7cm}
\end{minipage}
\begin{minipage}{0.48\linewidth}
\epsfig{file=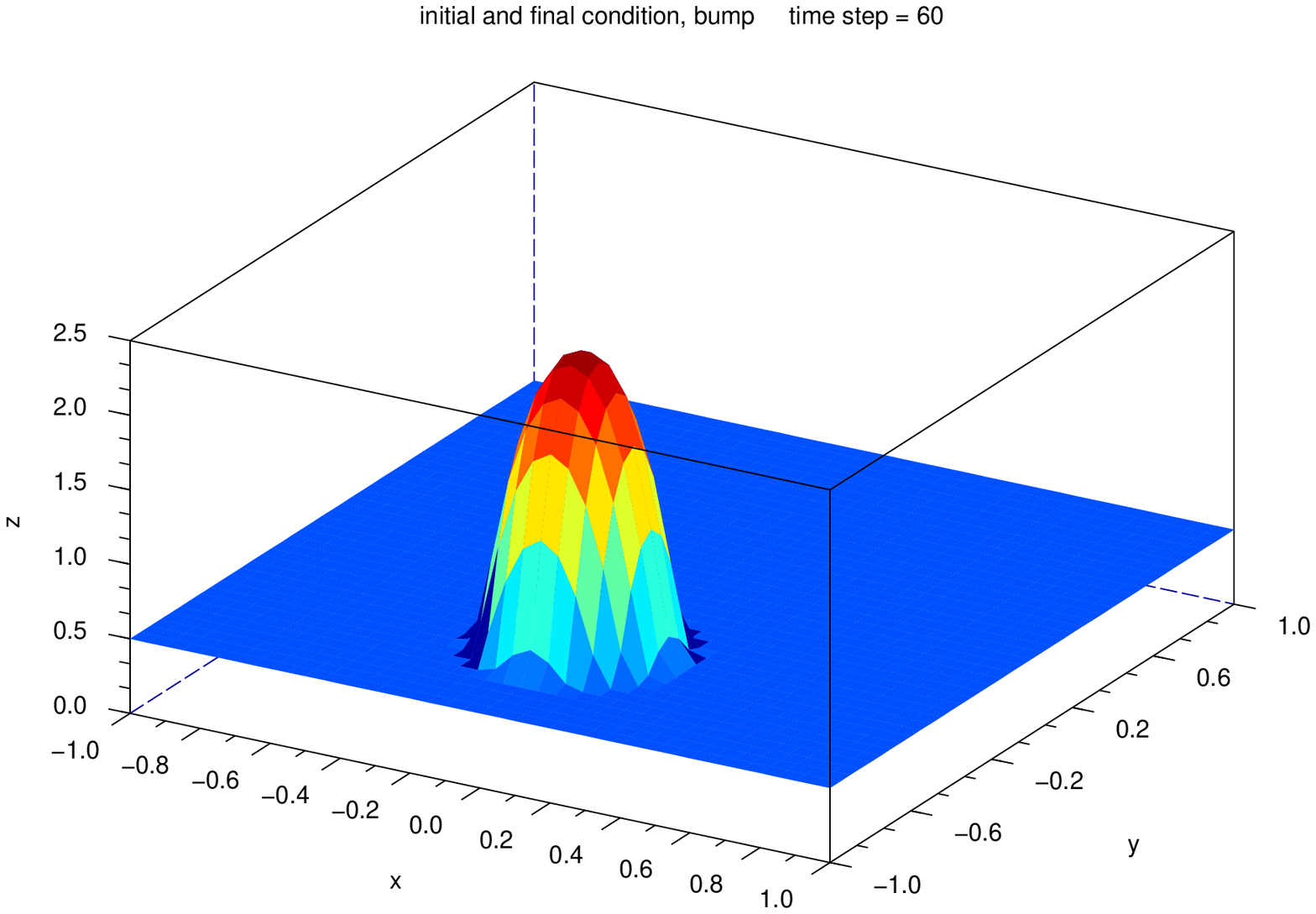,width=7cm}
\end{minipage}
\caption{Bump transport for $\beta=0.5$}
\label{beta05}
\end{center}
\end{figure}

\begin{figure}[H]
\begin{center}
\begin{minipage}{0.48\linewidth}
\epsfig{file=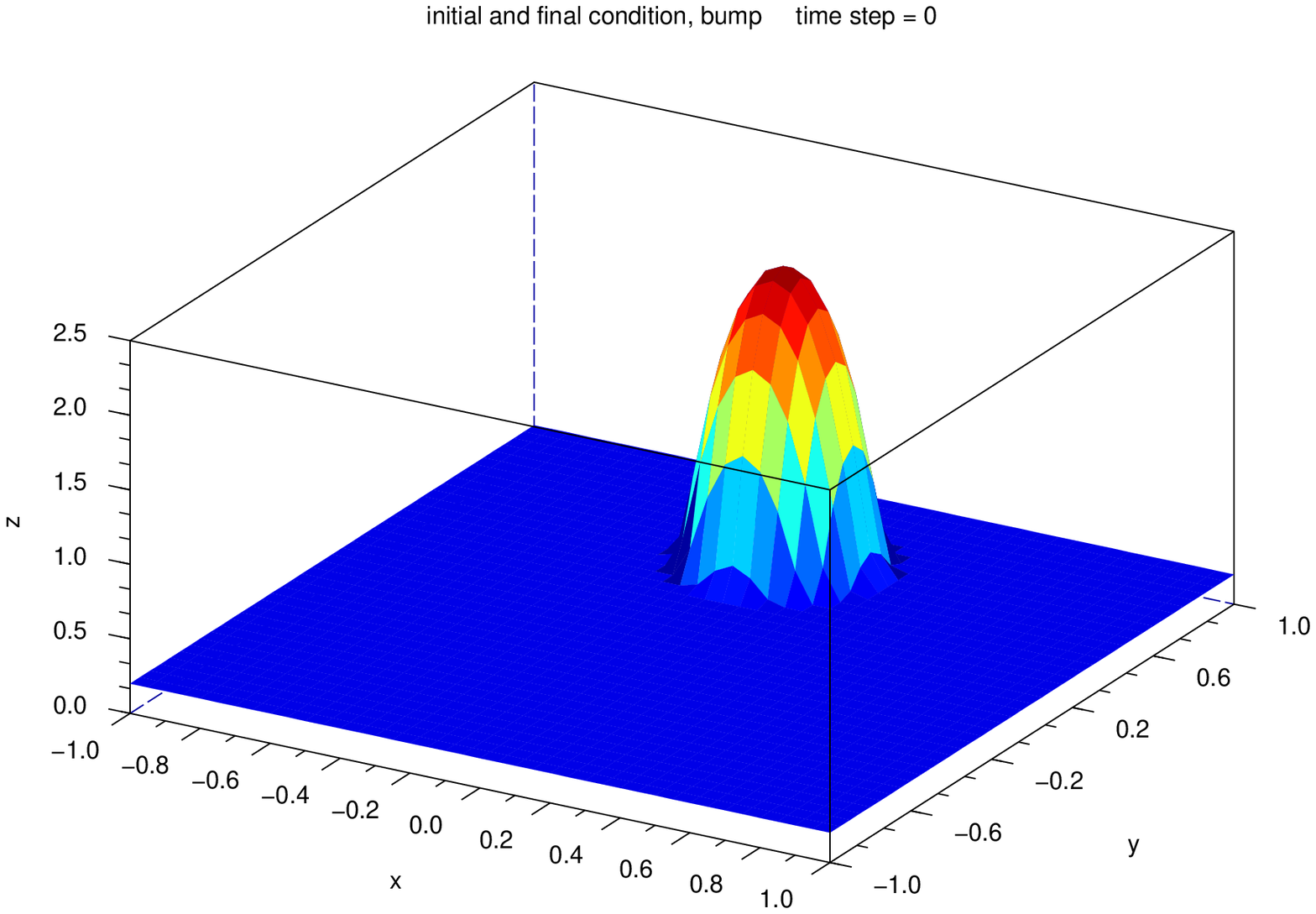,width=7cm}
\end{minipage}
\begin{minipage}{0.48\linewidth}
\epsfig{file=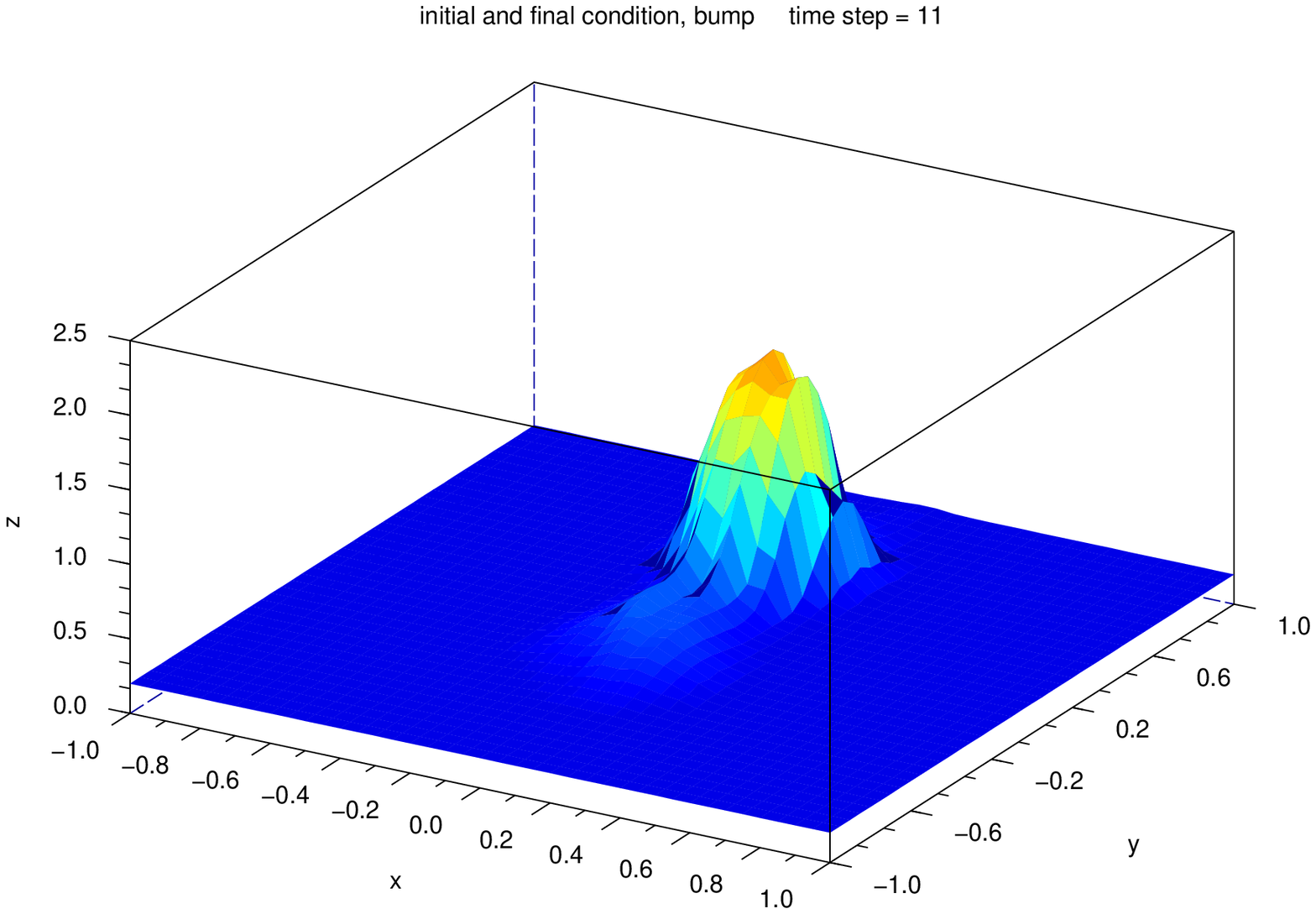,width=7cm}
\end{minipage}
\begin{minipage}{0.48\linewidth}
\epsfig{file=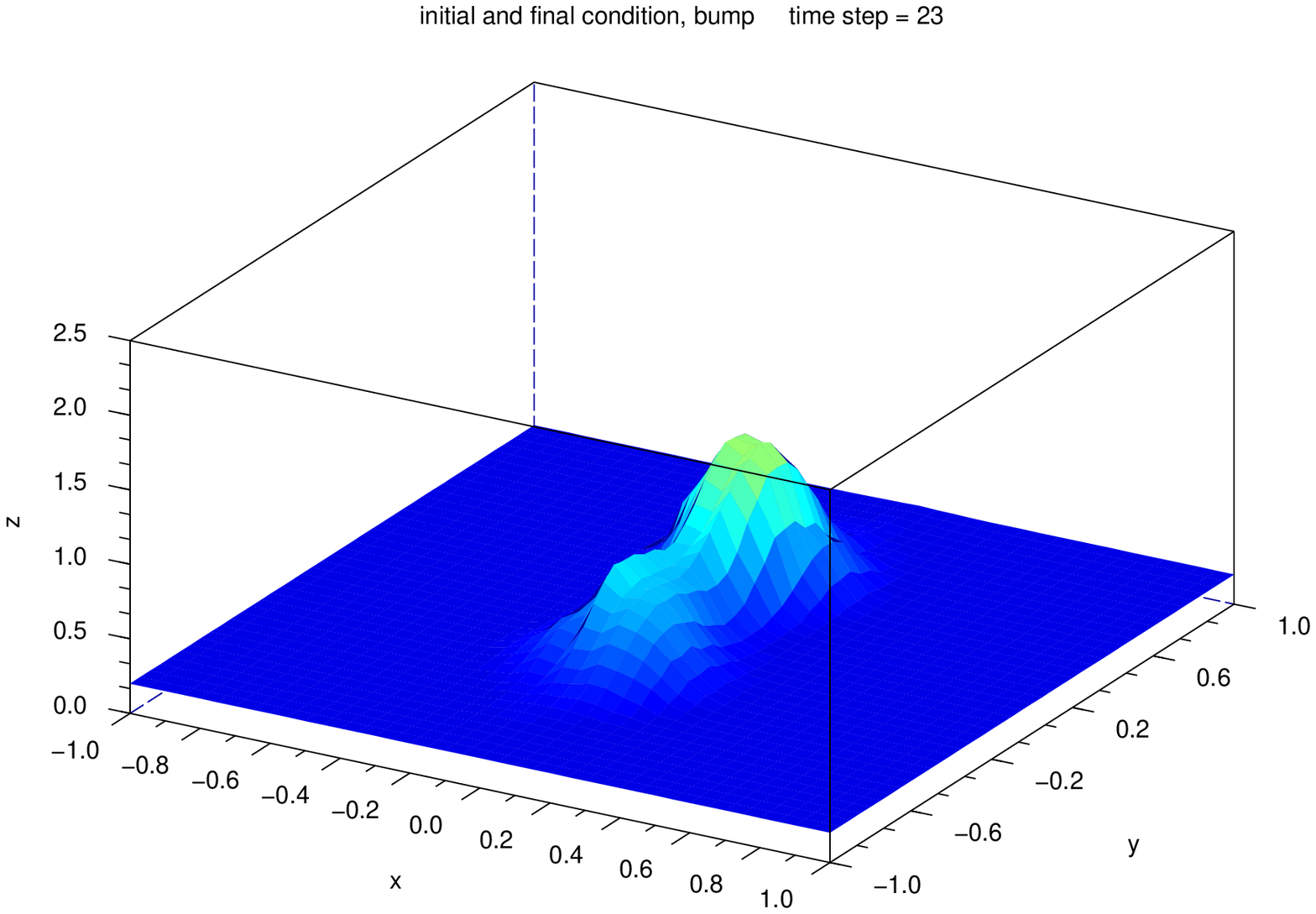,width=7cm}
\end{minipage}
\begin{minipage}{0.48\linewidth}
\epsfig{file=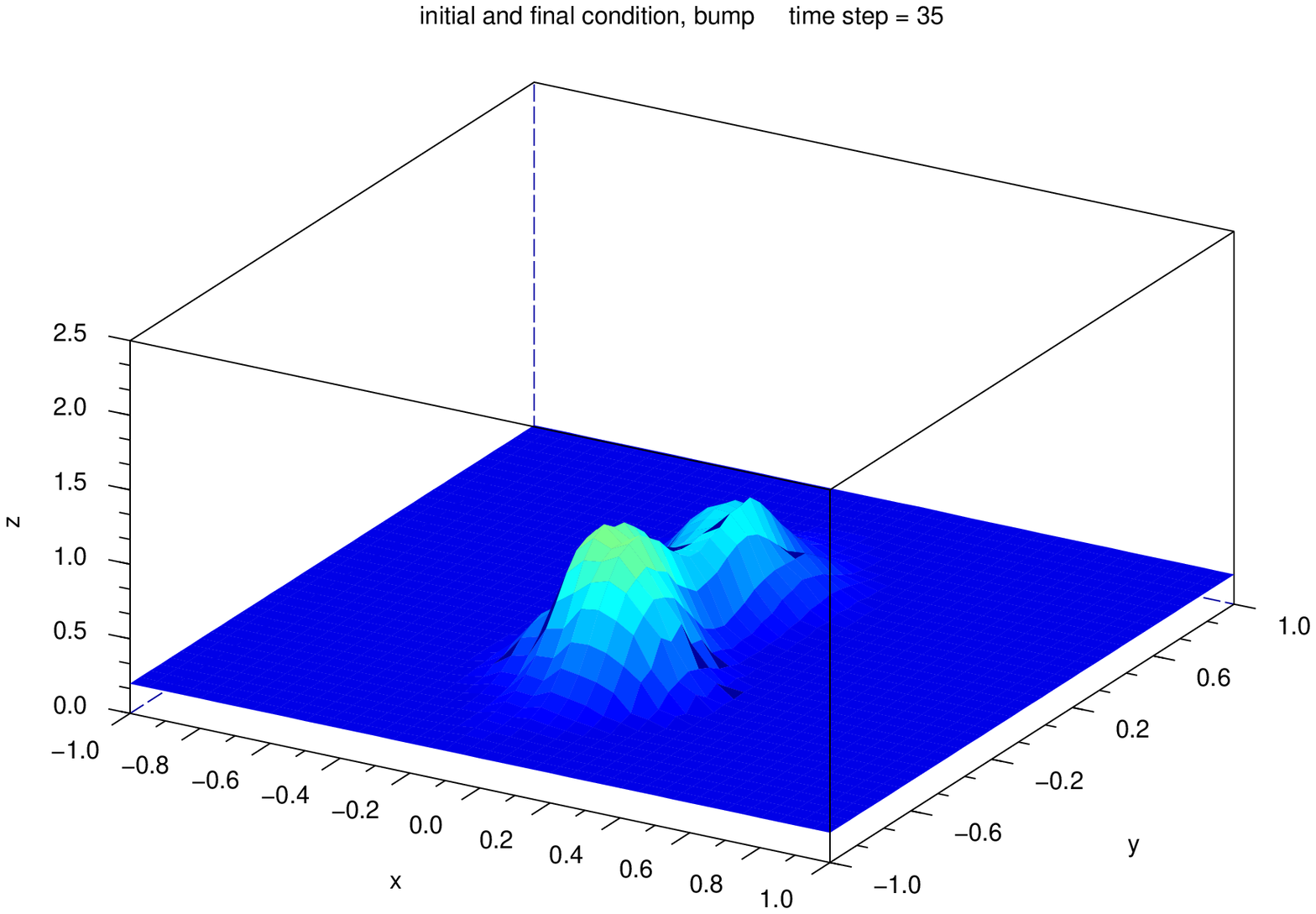,width=7cm}
\end{minipage}
\begin{minipage}{0.48\linewidth}
\epsfig{file=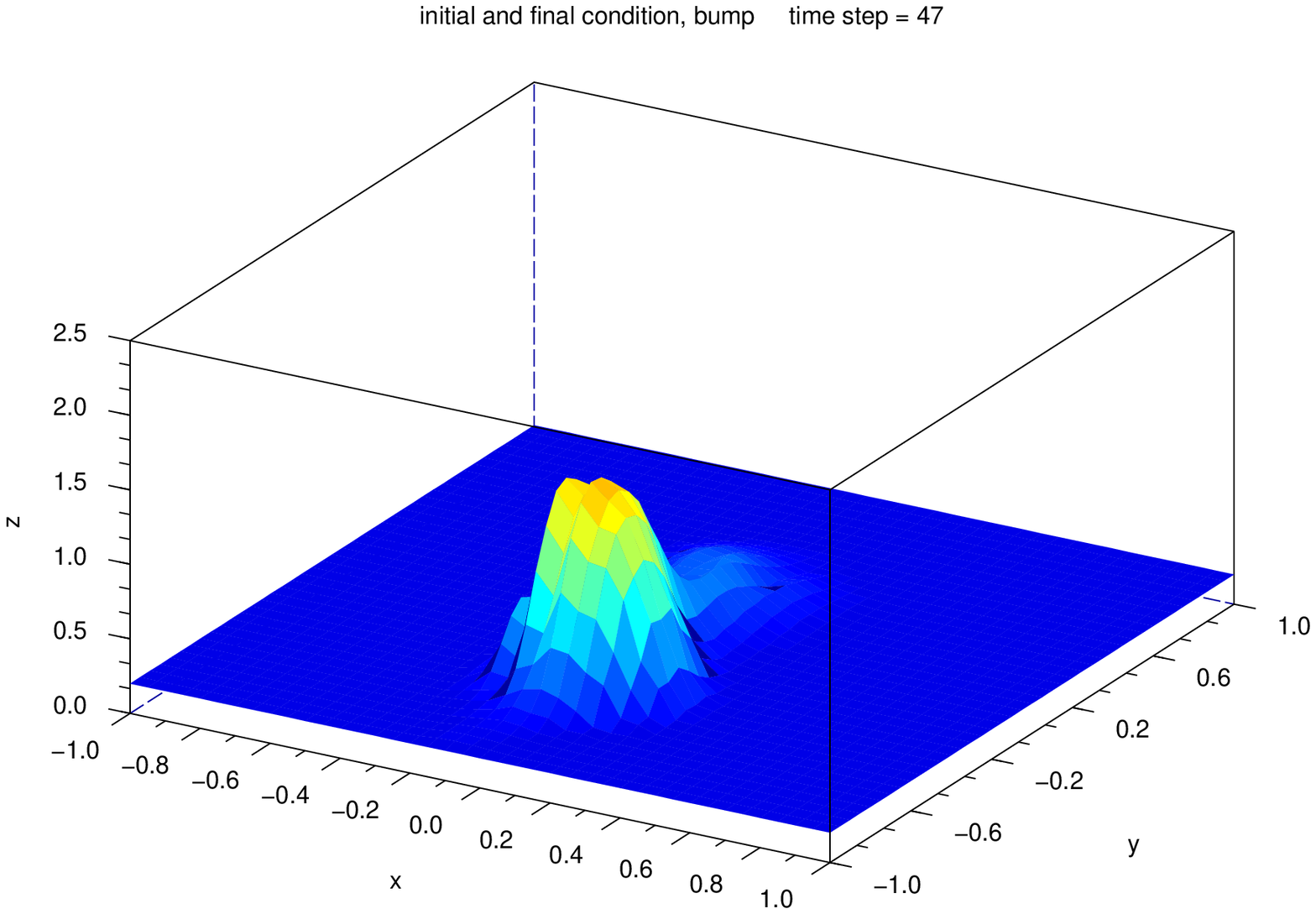,width=7cm}
\end{minipage}
\begin{minipage}{0.48\linewidth}
\epsfig{file=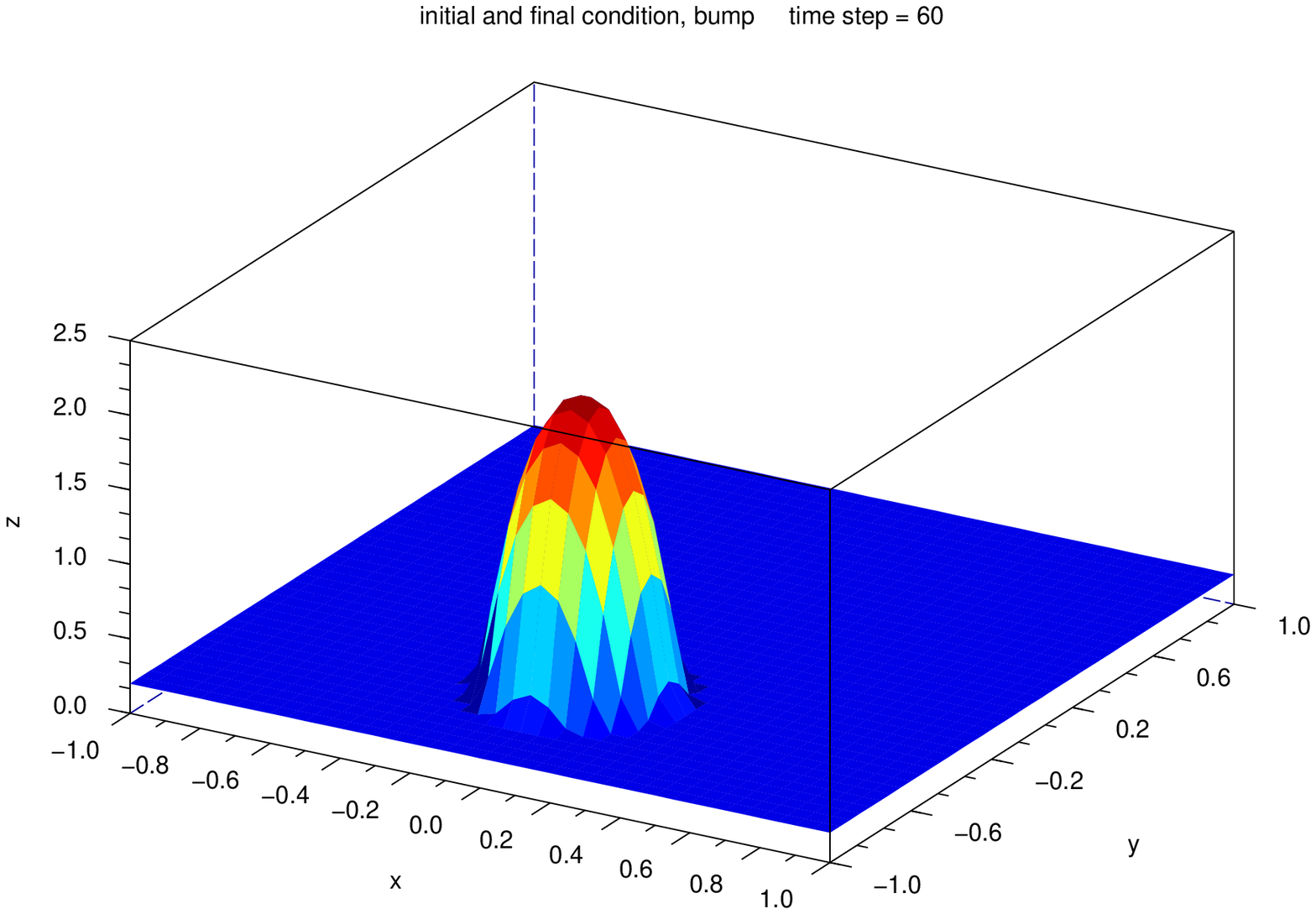,width=7cm}
\end{minipage}
\caption{Bump transport for $\beta=0.2$}
\label{beta02}
\end{center}
\end{figure}

\begin{figure}[H]
\begin{center}
\begin{minipage}{0.48\linewidth}
\epsfig{file=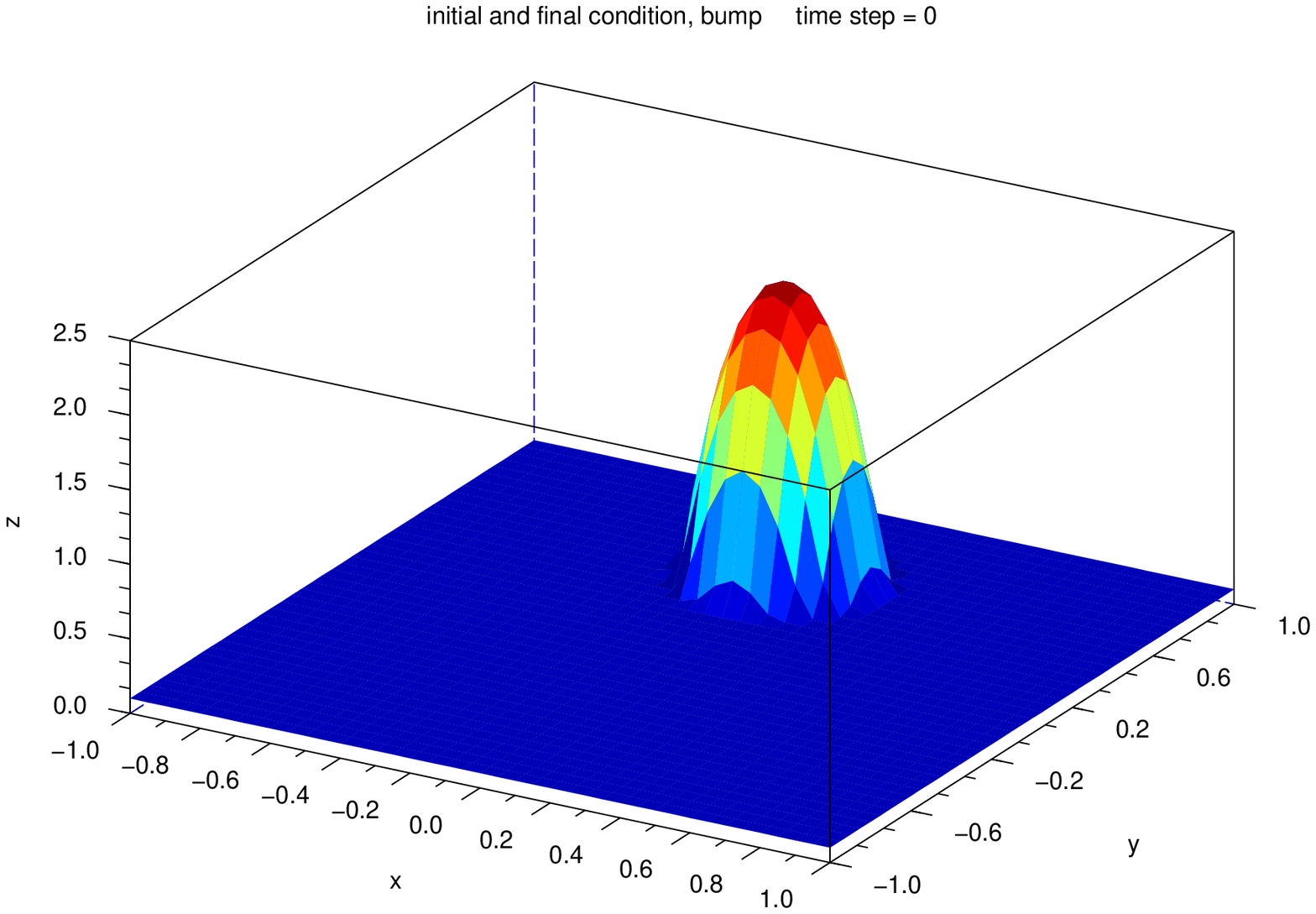,width=7cm}
\end{minipage}
\begin{minipage}{0.48\linewidth}
\epsfig{file=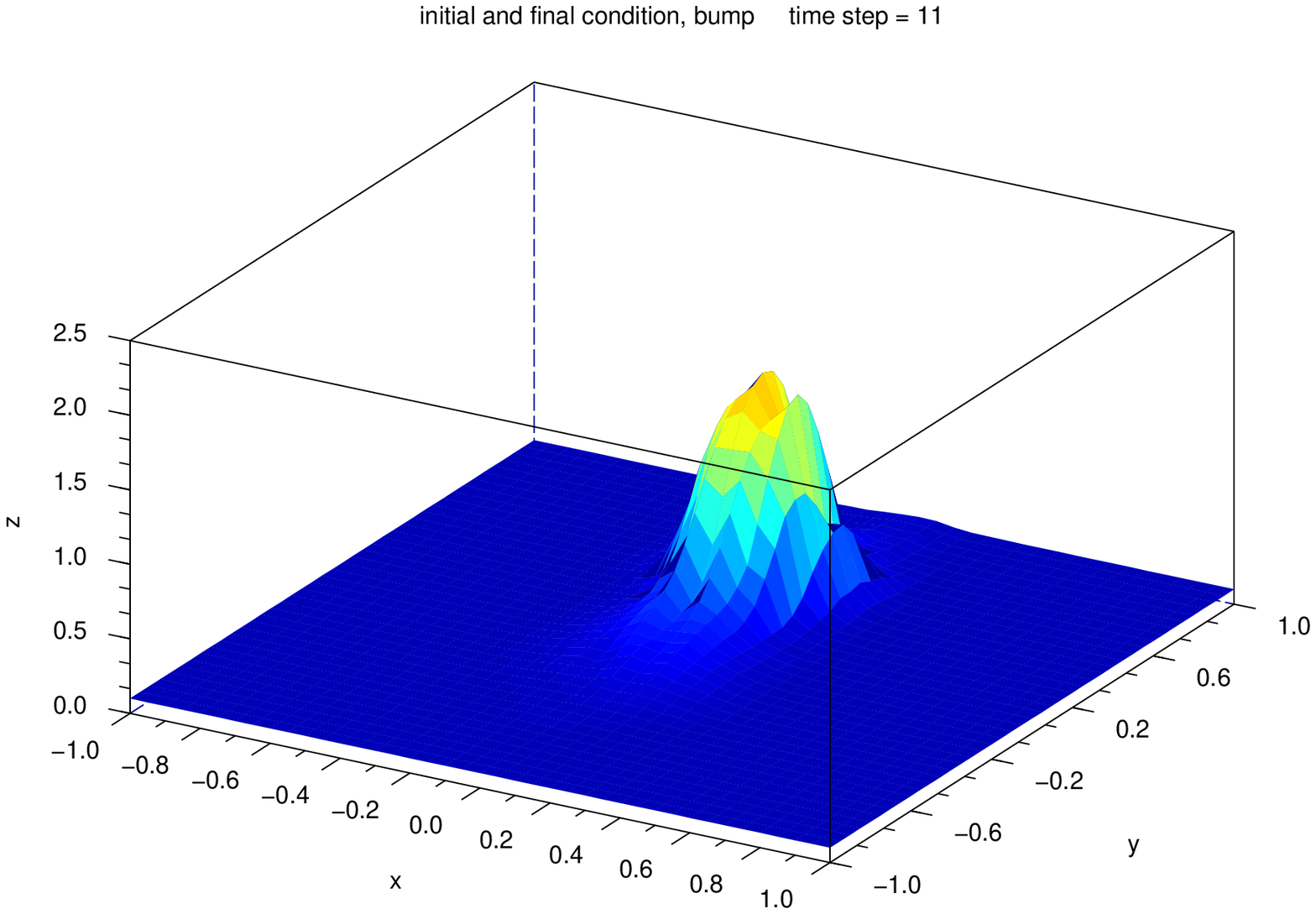,width=7cm}
\end{minipage}
\begin{minipage}{0.48\linewidth}
\epsfig{file=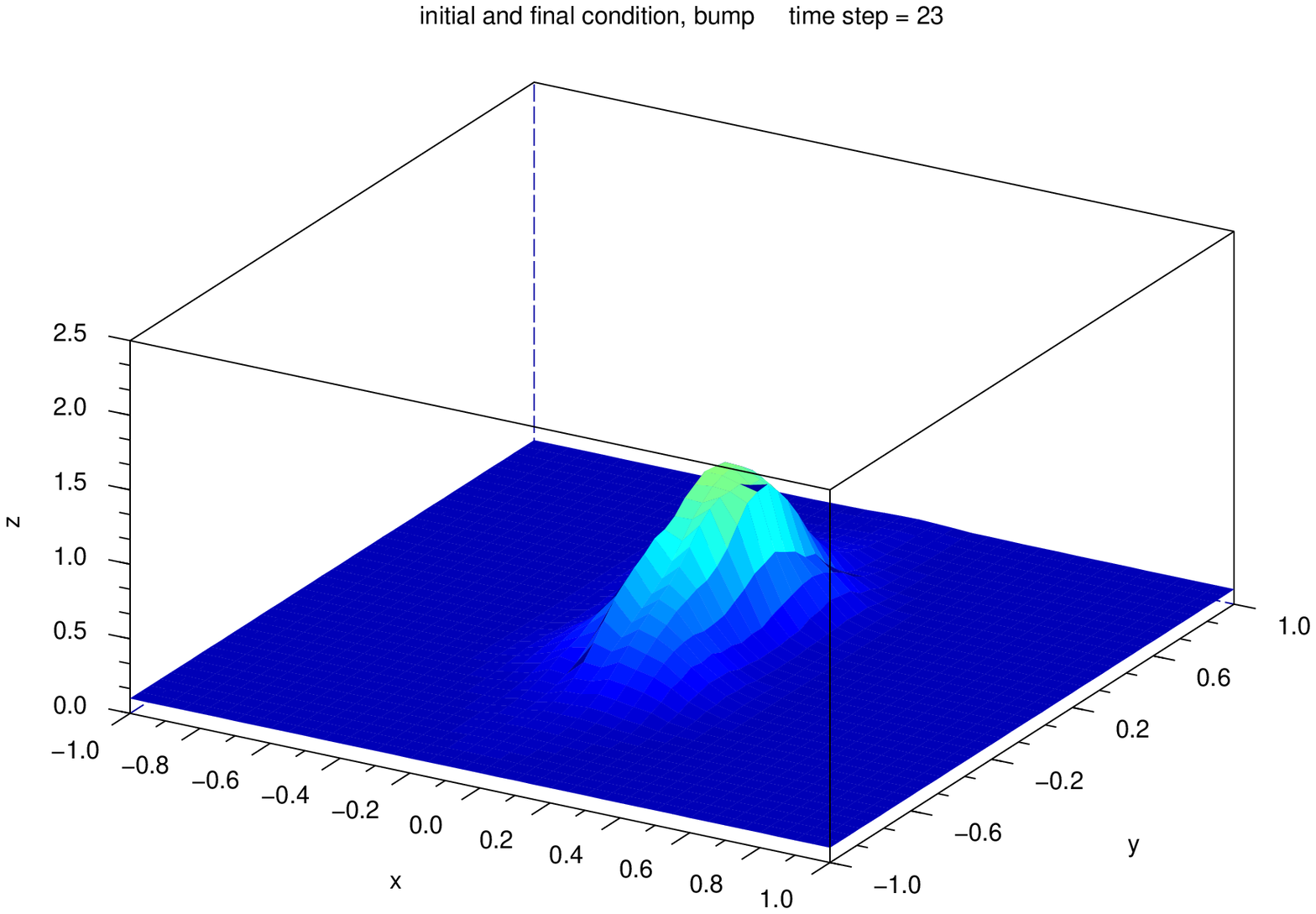,width=7cm}
\end{minipage}
\begin{minipage}{0.48\linewidth}
\epsfig{file=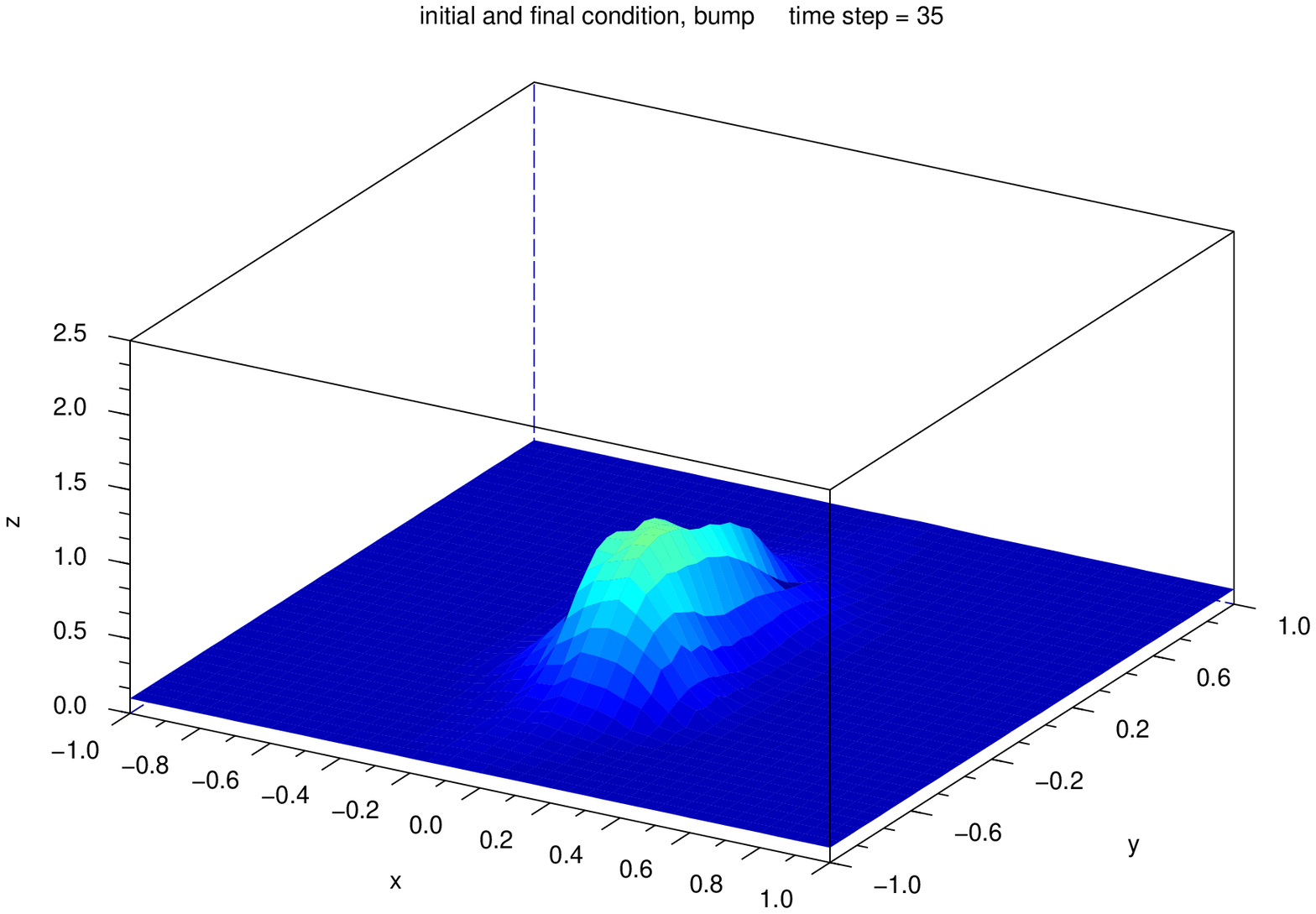,width=7cm}
\end{minipage}
\begin{minipage}{0.48\linewidth}
\epsfig{file=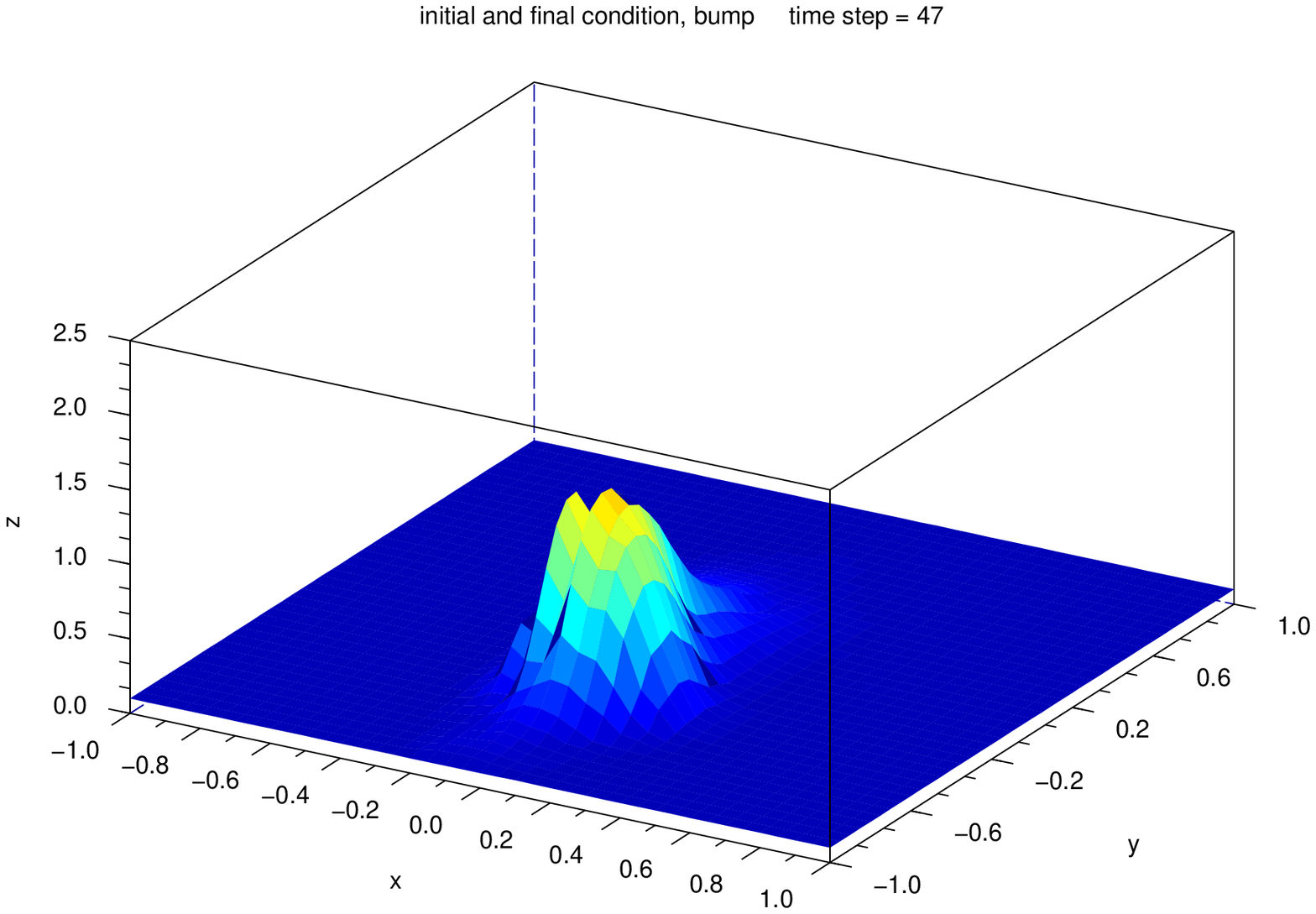,width=7cm}
\end{minipage}
\begin{minipage}{0.48\linewidth}
\epsfig{file=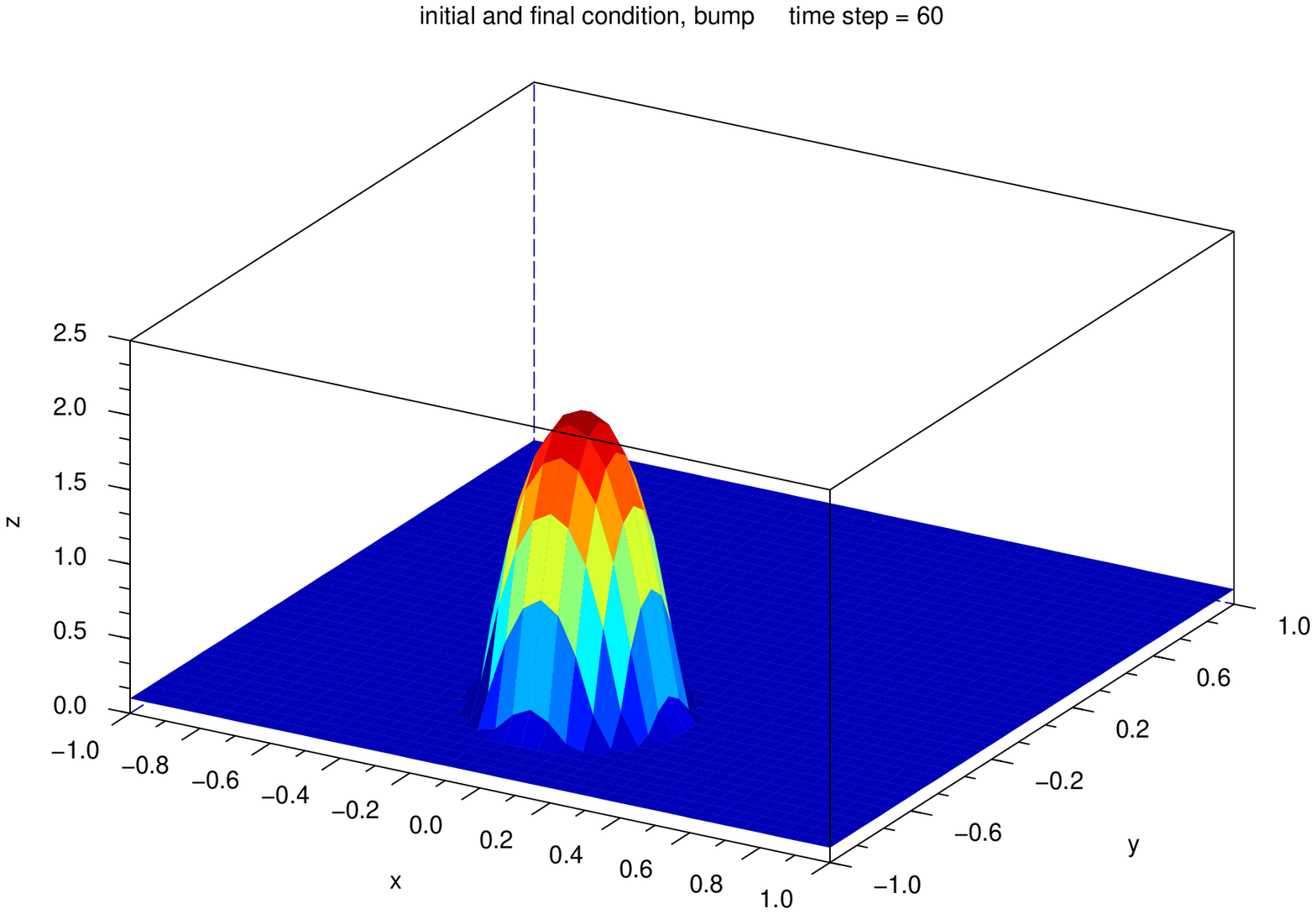,width=7cm}
\end{minipage}
\caption{Bump transport for $\beta=0.1$}
\label{beta01}
\end{center}
\end{figure}

\begin{figure}[H]
\begin{center}
\begin{minipage}{0.48\linewidth}
\epsfig{file=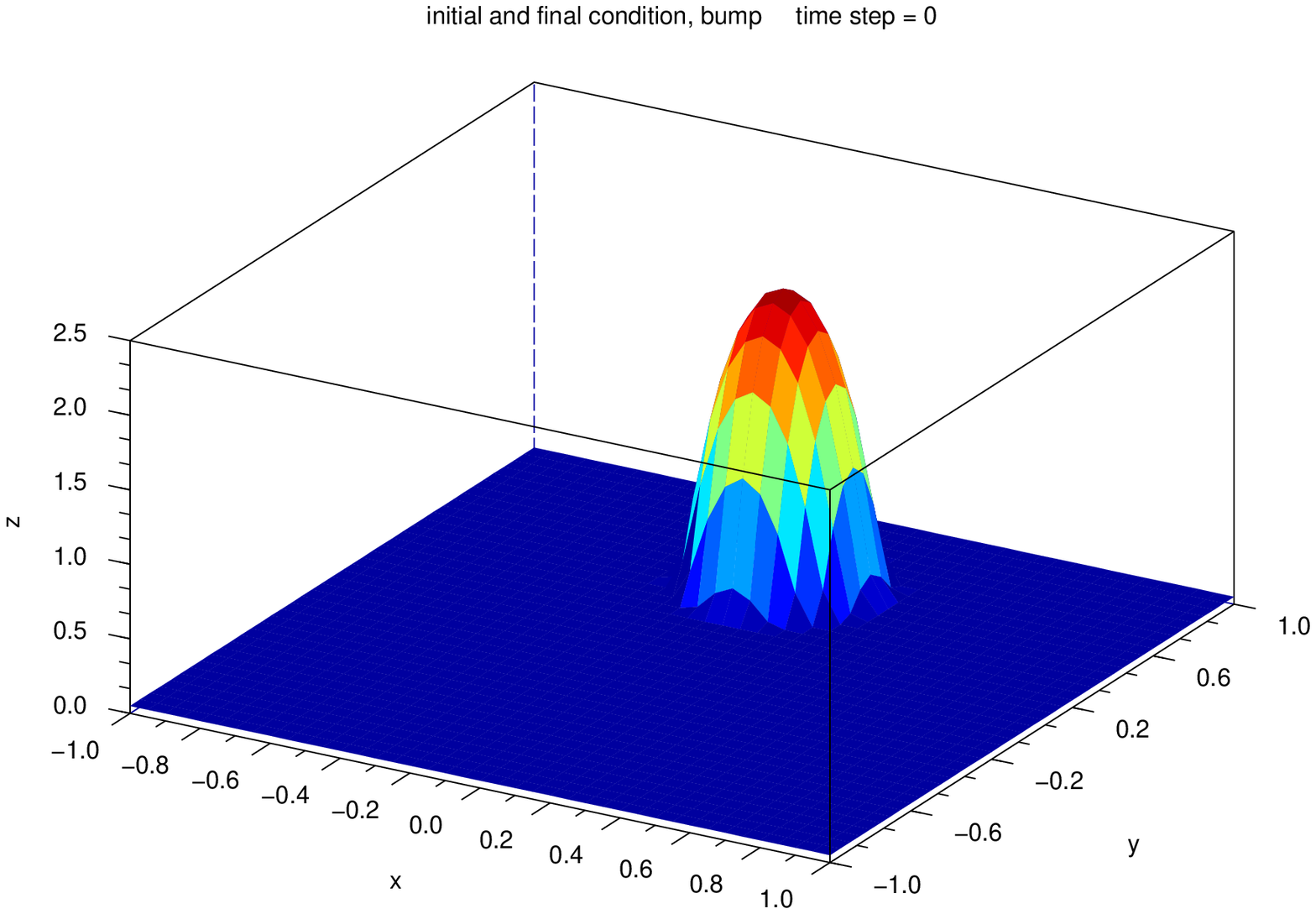,width=7cm}
\end{minipage}
\begin{minipage}{0.48\linewidth}
\epsfig{file=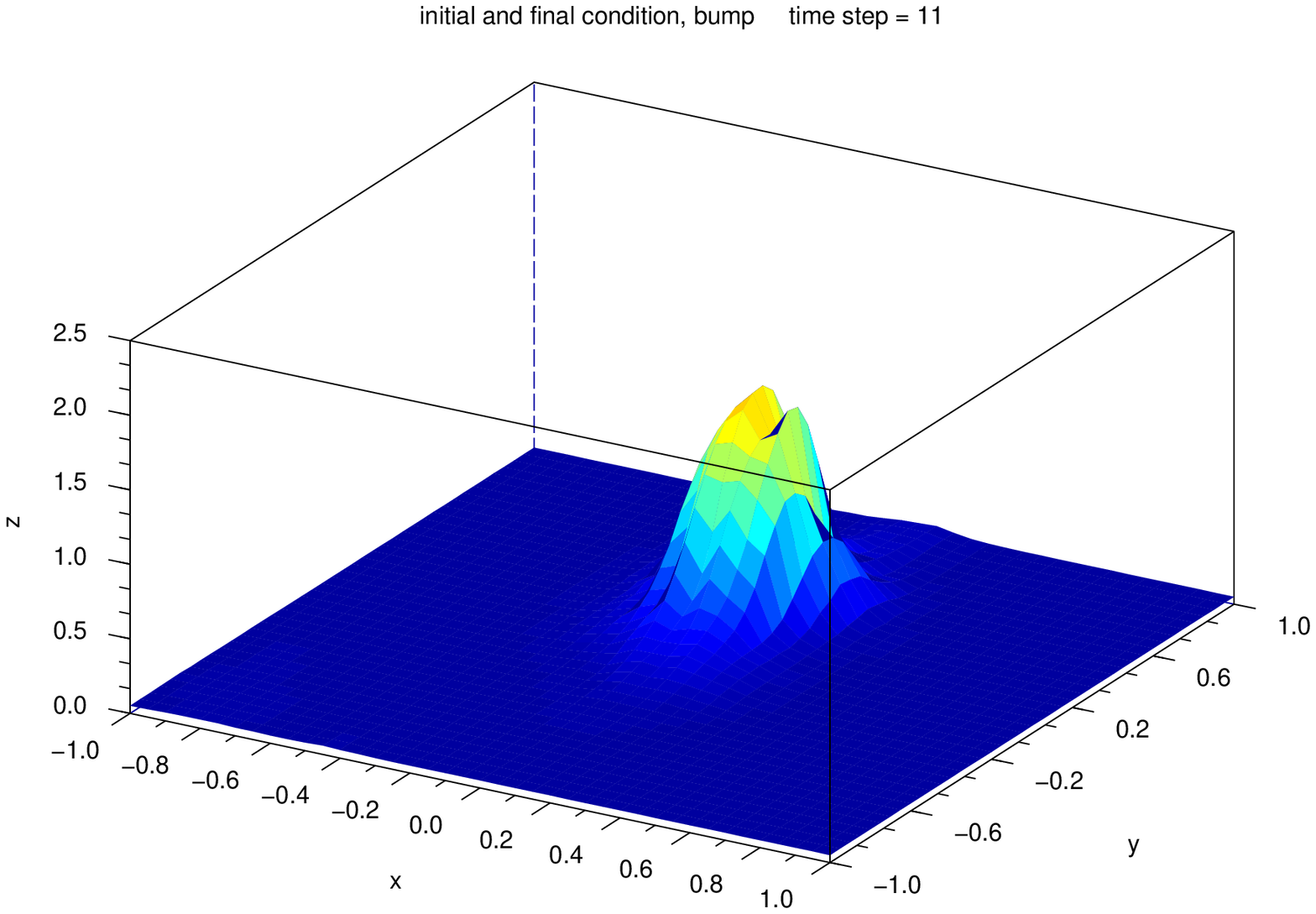,width=7cm}
\end{minipage}
\begin{minipage}{0.48\linewidth}
\epsfig{file=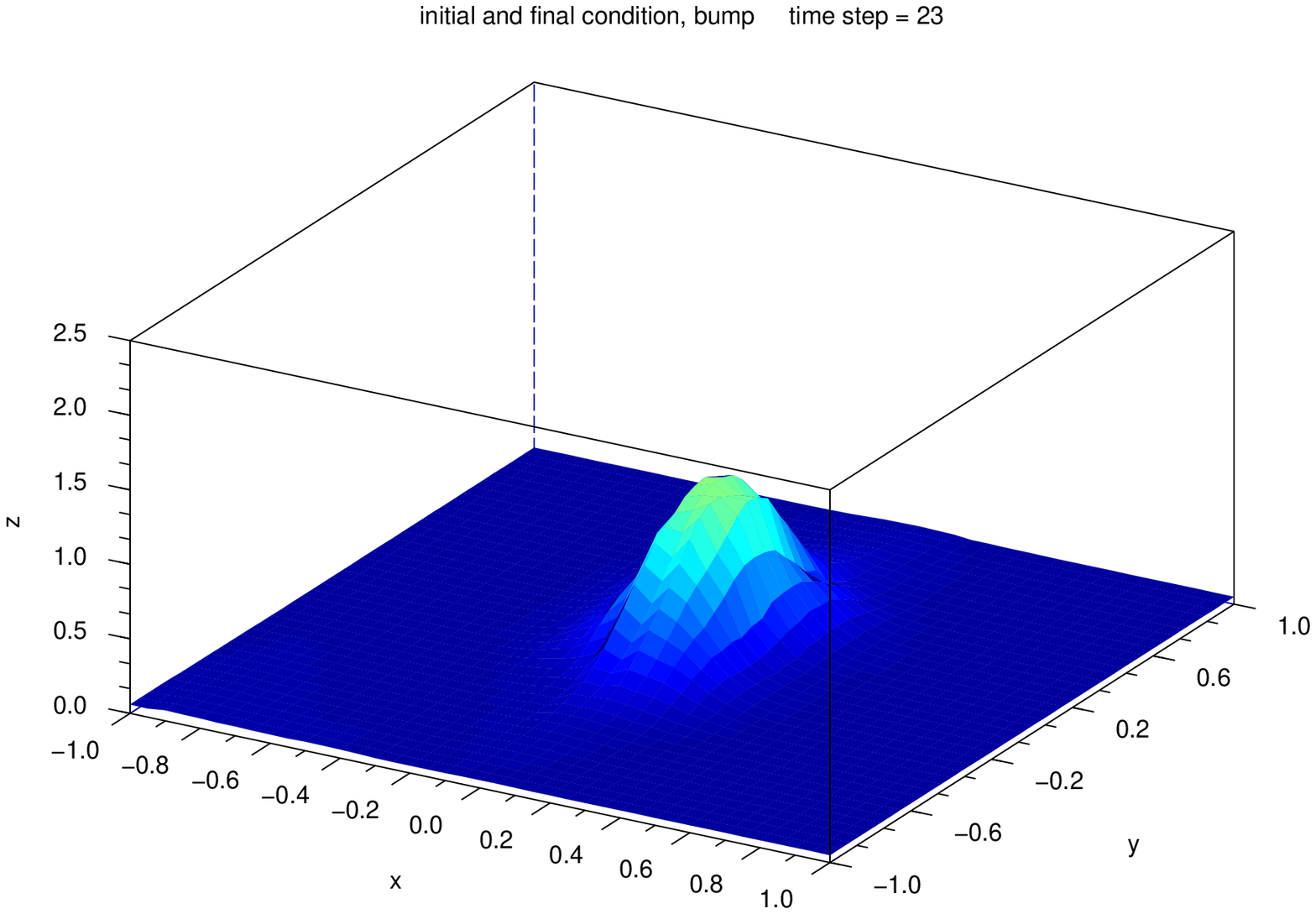,width=7cm}
\end{minipage}
\begin{minipage}{0.48\linewidth}
\epsfig{file=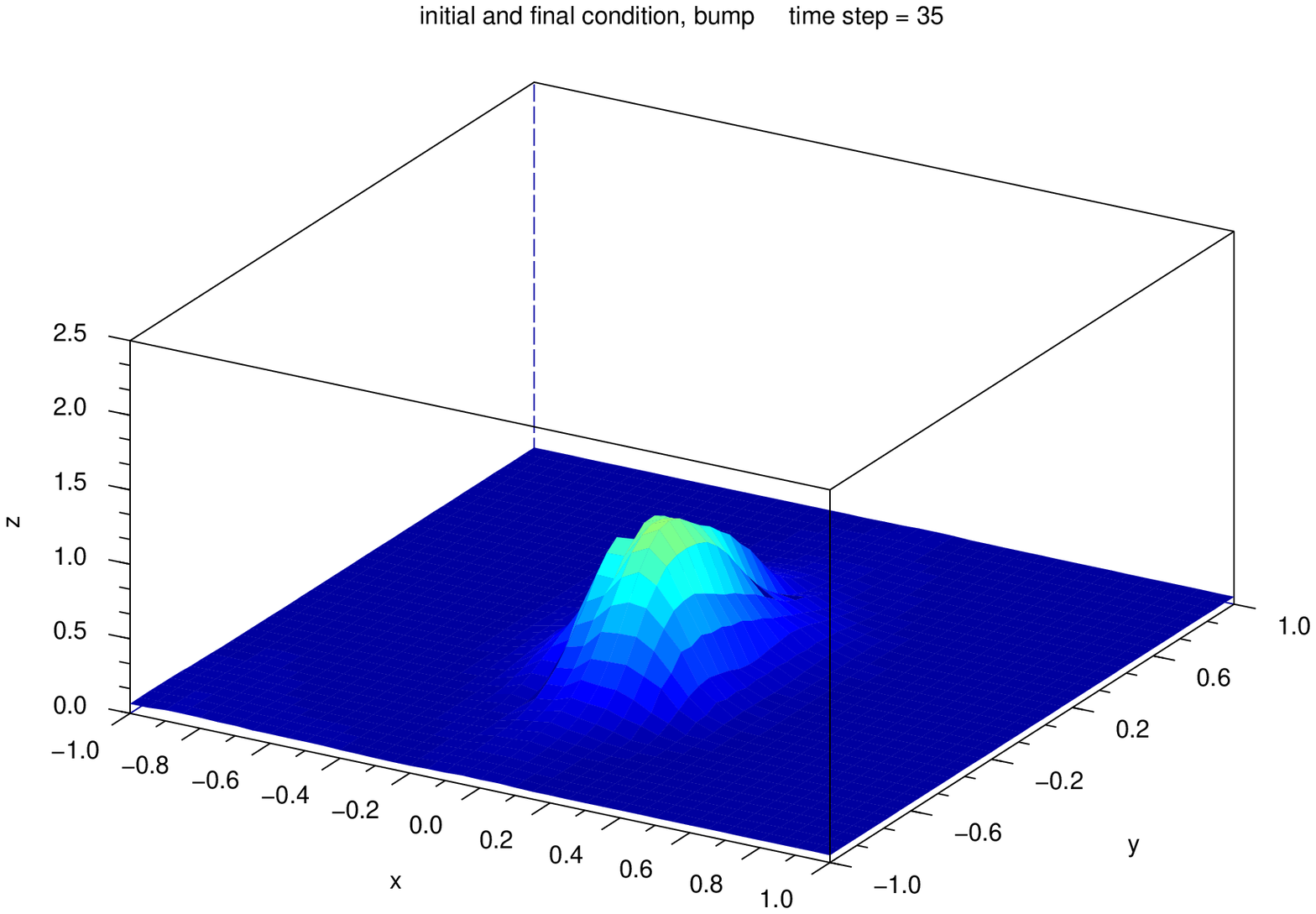,width=7cm}
\end{minipage}
\begin{minipage}{0.48\linewidth}
\epsfig{file=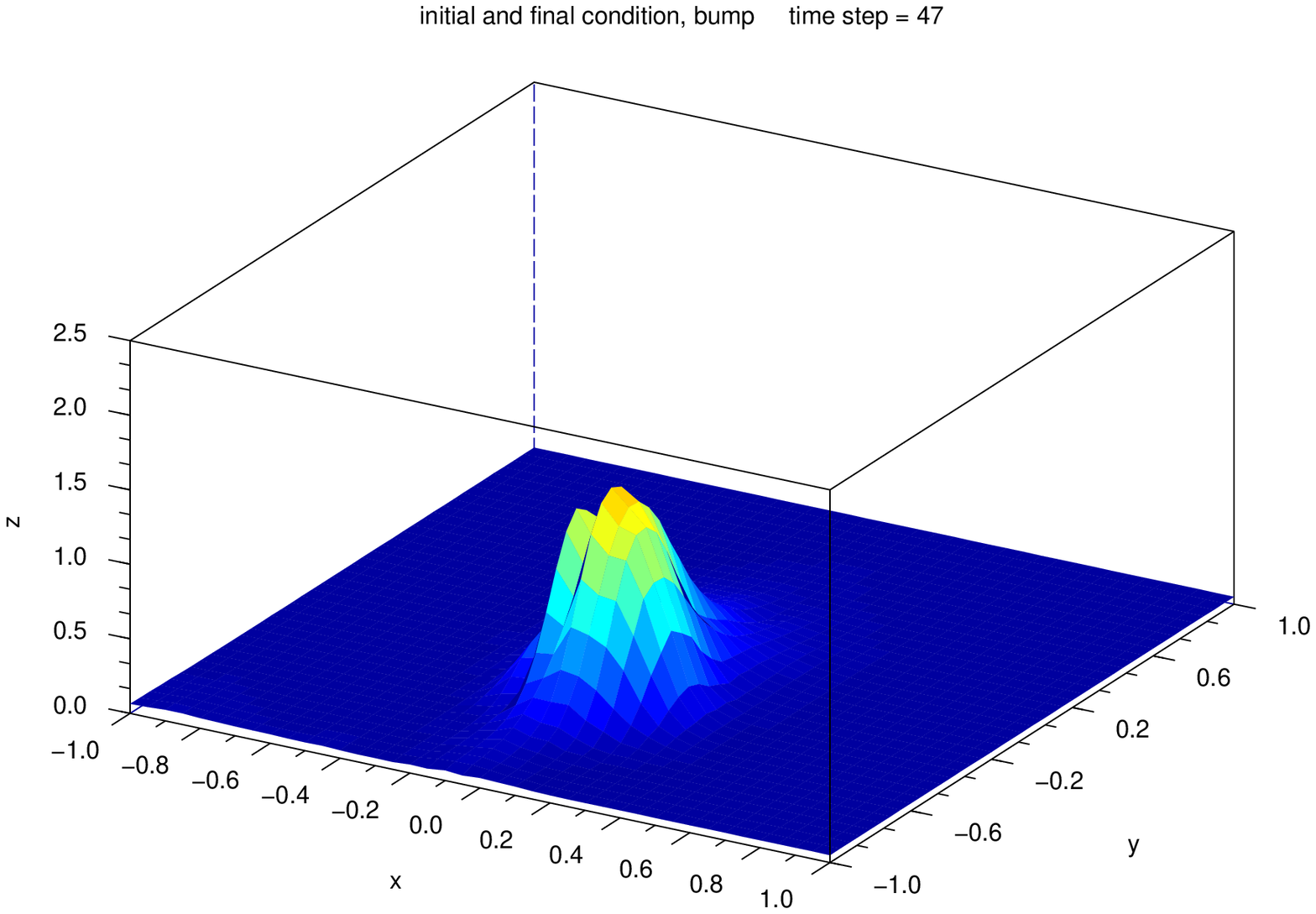,width=7cm}
\end{minipage}
\begin{minipage}{0.48\linewidth}
\epsfig{file=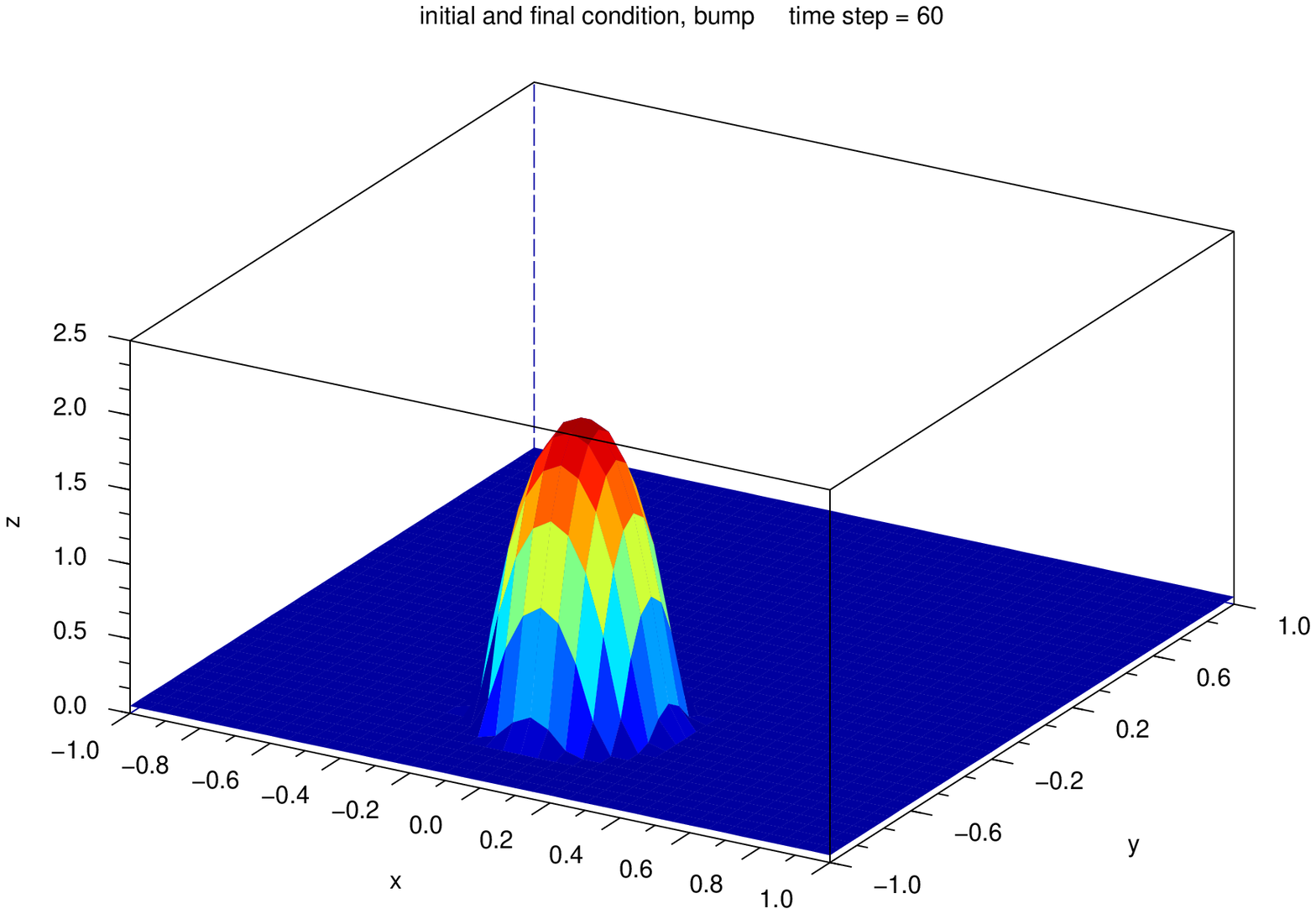,width=7cm}
\end{minipage}
\caption{Bump transport for $\beta=0.05$}
\label{beta005}
\end{center}
\end{figure}

In figures \ref{beta1}-\ref{beta005}, the shape of the bump is presented at the time-steps $0$,
$12$,
$24$, $36$, $48$, and $60$, the first one is the initial shape, and the last one is the final shape.
The next figure \ref{beta1} describe the bump transport for $\beta = 1$. Remark that in this case
this transport is very similar to the Dacorogna-Moser one.
Then figures \ref{beta05}, \ref{beta02}, \ref{beta01}, and \ref{beta005} describe the bump transport
for $\beta = 0.5$, $0.2$, $0.1$, and $0.05$.
Remark that when $\beta$ is small, the linear system becomes very ill-conditioned. Indeed for
$\beta \leq 0.1$ usual preconditioners like the classical IC0 one are useless. In this case the
parallel Gram-Schmidt least squares preconditioner (DIAG + LS CGS OPT)
developed in \cite{julien_these,julien_precond,julien_precond_paral} is used. Then the linear
system is solved in parallel using $48$ processors. The global tolerance for the
iterative scheme
developed in section \ref{algo} is set to $0.01$.

The number of iterations, and the CPU time for each value of $\beta$ is given in table
\ref{tableau}.
For $\beta < 0.05$ the algorithm did not converge.
\begin{table}[H]
\begin{center}
\begin{tabular}{|c|c|c|c|}
\hline
$\beta$ & nb. iterations & CPU time [s] & residue \\
\hline
\hline
1       &    4         &     128      & $8.62 \, 10^{-4}$ \\
\hline
0.5     &    5         &     149      & $9.23 \, 10^{-3}$ \\
\hline
0.2     &    8         &     201      & $4.41 \, 10^{-3}$ \\
\hline
0.1     &    10        &     229      & $9.47 \, 10^{-3}$ \\
\hline
0.05    &    55        &     1019     & $7.77 \, 10^{-3}$ \\
\hline
\end{tabular}
\caption{CPU time and number of iterations}
\end{center}
\label{tableau}
\end{table}
\subsection{The left ventricle motion}
\label{sec:ventricule}
The iterative strategy  described in Section 2 is then used to compute an approximated
solution, and to
reconstruct the systole to diastole images of a slice of a left ventricle.
Ten time steps have been used to compute the solution, and 10000 degrees
of freedom for the time-space least squares finite element. The approximated
fixed point algorithm converges in about 10 iterations with an accuracy of about
$10^{-4}$. In this case the usual IC0 preconditioner is sufficient; this is essentially due to the
fact that there is no large region in the domain with a very low density $\rho$.
In the next figure \ref{systole}, the initial image and the final image are presented.
\begin{figure}[H]
\begin{center}
\begin{minipage}{0.48\linewidth}
\epsfig{file=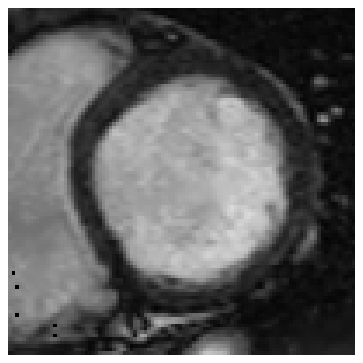, width=5cm}
\end{minipage}
\begin{minipage}{0.48\linewidth}
\epsfig{file=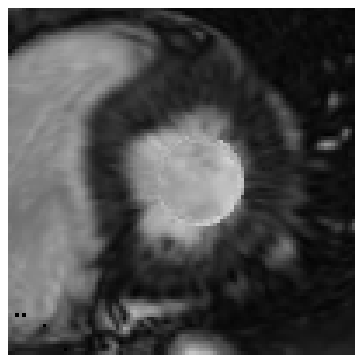,width=5cm}
\end{minipage}
\caption{End of diastole of a left ventricular (a), of systole (b) }\label{systole}
\end{center}
\end{figure}
In the following  figure \ref{syst_5}, two intermediate times $1/3$ and $2/3$ are shown.
\begin{figure}[H]
\begin{center}
\begin{minipage}{0.48\linewidth}
\epsfig{file=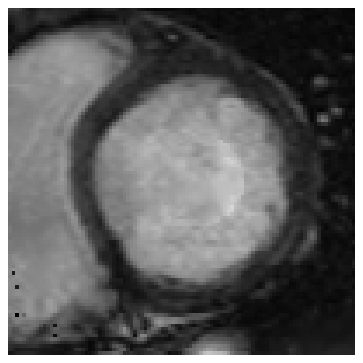, width=5cm}
\end{minipage}
\begin{minipage}{0.48\linewidth}
\epsfig{file=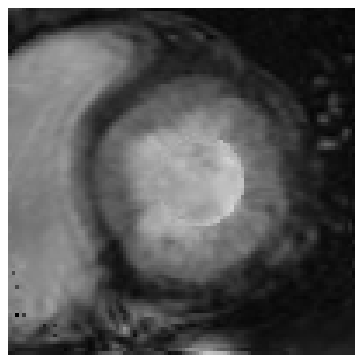,width=5cm}
\end{minipage}
\caption{Time step 3 and 6}\label{syst_5}
\end{center}
\end{figure}
To summarize, in this work, we present a fixed point algorithm
for the computation of the time dependent optimal mass transportation problem, allowing to  handle
the images tracking
problem. The efficiency of the method has been tested with some 2D examples.
%
%

\end{document}